\documentclass[a4paper,fleqn]{cas-sc}

\usepackage[numbers]{natbib}

\usepackage{amsmath}
\usepackage{graphicx}
\usepackage{amsthm}
\usepackage{bm}
\usepackage{enumitem}
\usepackage{soul,color}
\usepackage{xcolor}
\usepackage{mathtools}

\def\P{{\cal P}}

\def\E{{\mathbb E}}
\def\O{{\cal O}}
\def\R{{\mathbb R}}
\def\N{{\cal N}}
\def\C{{\cal C}}

\newdefinition{definition}{Definition}
\newdefinition{remark}{Remark}
\newdefinition{example}{Example}
\newtheorem{theorem}{Theorem}
\newtheorem{lemma}{Lemma}
\newtheorem{proposition}{Proposition}
\newtheorem{corollary}{Corollary}

\newcommand{\hlnew}[1]{{#1}}

\newcommand{\norm}[1]{\|#1\|} 
\newcommand\abs[1]{\left|#1\right|}
\newcommand\numberthis{\addtocounter{equation}{1}\tag{\theequation}}

\DeclareMathOperator*{\vect}{vec}

\DeclareMathOperator*{\tr}{tr}
\DeclareMathOperator*{\diag}{diag}

\begin{document}

\ExplSyntaxOn
\keys_set:nn { stm / mktitle } { nologo }
\ExplSyntaxOff

\let\printorcid\relax
\let\WriteBookmarks\relax

\def\floatpagepagefraction{1}
\def\textpagefraction{.001}

\shorttitle{Perturbation expansions and error bounds for TSVD}
\shortauthors{Vu et. al.}

\title [mode = title]{Perturbation expansions and error bounds for the truncated singular value decomposition}                      

\tnotetext[1]{This work was partially supported by the National Science Foundation grant CCF-1254218.}


\author[1]{Trung Vu}
\cormark[1]
\ead{vutru@oregonstate.edu}
\address[1]{School of Electrical Engineering and Computer Science, Oregon State University, Corvallis, OR 97331-5501, USA}
\cortext[cor1]{Corresponding author}

\author[2]{Evgenia Chunikhina}
\ead{chunikhina@pacificu.edu}
\address[2]{Department of Mathematics and Computer Science, Pacific University, Forest Grove, OR 97116-1797, USA}

\author[1]{Raviv Raich}
\ead{raich@eecs.oregonstate.edu}

\begin{abstract}
Truncated singular value decomposition is a reduced version of the singular value decomposition in which only a few largest singular values are retained. This paper presents a \hlnew{novel} perturbation analysis for the truncated singular value decomposition for real matrices. First, we \hlnew{describe} perturbation expansions for the singular value truncation of order $r$. We extend perturbation results for the singular subspace decomposition to derive the first-order perturbation expansion of the truncated operator about a matrix with rank greater than or equal to $r$. Observing that the first-order expansion can be greatly simplified when the matrix has exact rank $r$, we further show that the singular value truncation admits a simple second-order perturbation expansion about a rank-$r$ matrix. Second, we introduce the first-known error bound on the linear approximation of the truncated singular value decomposition of a perturbed rank-$r$ matrix. Our bound only depends on the least singular value of the unperturbed matrix and the norm of the perturbation matrix. Intriguingly, while the singular subspaces are known to be extremely sensitive to additive noises, the \hlnew{newly established} error bound holds universally for perturbations with arbitrary magnitude. \hlnew{Finally, we demonstrate an application of our results to the analysis of the mean squared error associated with the TSVD-based matrix denoising solution.}
\end{abstract}




\maketitle

\section{Introduction}

The singular value decomposition (SVD) is an invaluable tool for matrix analysis and the truncated singular value decomposition (TSVD) offers a formal approach for a rank-restricted optimal approximation of matrices by replacing the smallest singular values by zeros in the SVD of a matrix. TSVD has numerous applications in science, engineering, and math with examples including linear system identification \cite{markovsky2008structured,markovsky2005application}, collaborative filtering \cite{candes2009exact,jain2010guaranteed}, low-rank matrix denoising \cite{shabalin2013reconstruction,yang2016rate}, data compression \cite{wei2001ecg}, and numerical partial differential equations \cite{li2011ill}. In addition, TSVD is well-known for solving classical discrete ill-posed problems \cite{hansen1987truncatedsvd, hansen1990truncated}. This paper is concerned with the effects of errors on the truncated singular value decomposition of a matrix.

Perturbation theory for the SVD studies the effect of variation in matrix entries on the singular values and the singular vectors of a matrix. 
Using perturbation bounds or perturbation expansions, one can characterize the difference between the SVD-related quantities associated with the perturbed matrix and those of the original matrix. 
The first perturbation bound on singular values was given by Weyl \cite{weyl1912asymptotische} in 1912, stating that no singular value can \hlnew{be} changed by more than the spectral norm of the perturbation. 
Later, Mirsky \cite{mirsky1960symmetric} showed that Weyl's inequality also holds for any unitarily-invariant norm. \hlnew{Perturbation bounds for singular vectors are often established in the context of singular subspace decomposition. In 1970, Davis and Kahan} \cite{davis1970rotation} \hlnew{introduced a fundamental bound on the distance between the subspaces spanned by a group of eigenvectors and their perturbed versions based the ratio between the perturbation level and the eigengap. This result is also referred as the so-called $\sin \Theta$ theorem for symmetric matrices in the literature.} Shortly afterwards, Wedin \cite{wedin1972perturbation} generalized part of this result to cover non-symmetric matrices using the singular value decomposition, bounding changes in the left and right singular subspaces in terms of the singular value gap and the perturbation magnitude. In a recent work, Cai and Zhang \cite{cai2018rate} \hlnew{further established separate matching upper and lower bounds for the left and right singular subspaces.}
When the structure of the error is concerned, one may draw interest in perturbation expansions to approximate the perturbed quantity as a function of the perturbation matrix. As the perturbation decreases towards zero, the approximation is more accurate since the higher-order terms in the expansion become successively smaller. In 1973, Stewart \cite{stewart1973error} showed that there exists explicit expression of the perturbed subspaces in the bases of the unperturbed subspaces, which can be leveraged to obtain error bounds for certain characteristic subspaces associated with the SVD. 
This breakthrough result has started a long line of research on perturbation expansions and error bounds for the SVD, including the work of Stewart \cite{stewart1984second}, Sun \cite{ji1988note}, Li~\textit{et~al.} \cite{li1991unified}, Vaccaro \cite{vaccaro1994second}, Xu \cite{xu2002perturbation}, Liu~\textit{et~al.} \cite{liu2008first}, and more recently, Gratton~\textit{et~al.} \cite{gratton2016second}. 
Specifically, in \cite{stewart1984second}, \hlnew{Stewart utilized the bounding technique in} \cite{stewart1973error} \hlnew{and obtained a second-order perturbation expansion for the square of the smallest singular value of a matrix. In a different approach based on the theory of implicit functions, Sun} \cite{ji1988note} \hlnew{provided the first analytical expression for the second-order perturbation expansion of simple non-zero singular values of a matrix. One of the first significant results on perturbation expansion of singular subspaces was introduced by Li and Vaccaro in 1991. In }\cite{li1991unified}, \hlnew{the two authors analyzed a variety of subspace-based algorithms in array signal processing and developed the first-order perturbation expansion for the signal and orthogonal subspaces of the rank-deficient data matrix. Later on, Vaccaro} \cite{vaccaro1994second} \hlnew{extended this result to the second-order perturbation expansion of these subspaces. A more fine-grained analysis of the perturbation expansion for the individual singular vectors rather than the singular subspaces was given by Liu~\textit{et~al.}} \cite{liu2008first}, \hlnew{uncovering the fact that the signal subspace has an impact on the first-order approximation of the individual singular vectors, but not on the first-order approximation of the signal subspace spanned by these vectors.
We note that the aforementioned results on perturbation analysis of singular subspaces make an assumption that the unperturbed matrix is rank-deficient, i.e., all singular values corresponding to one of the singular subspaces are zero. In 2002, Xu} \cite{xu2002perturbation} \hlnew{relaxed this constraint by only requiring those singular values to be equally small. Recently, Gratton and Tshimanga} \cite{gratton2016second} \hlnew{were able to eliminate this constraint completely, presenting the second-order perturbation expansion for singular subspaces with no restriction on their corresponding singular values.}\footnote{The only constraint is the singular-value separation between the two subspaces.} \hlnew{It is notable that the last result is developed directly from those by Stewart in} \cite{stewart1973error}.
\hlnew{A more comprehensive description of the aforementioned results is given in Section~3. Interested readers can also find in-depth surveys on matrix perturbation theory in} \cite{stewart1990perturbation,Stewart90} and references therein.

\hlnew{The aforementioned results on perturbation analysis of the SVD is the fulcrum for the perturbation analysis of the TSVD. While the former characterizes the effect of perturbation on the singular values/singular subspaces of a matrix, the later studies the combined effect (from both singular values and singular subspaces) on the resulting reduced-rank matrix. Analyzing such an effect helps understand the local behavior of algorithms that utilize the low-rank optimal approximation of matrices, such as SVD-based channel estimation methods in multi-input multi-output (MIMO) systems} \cite{lindskog1999reduced,nicoli2005reduced,jing2012ml} \hlnew{and iterative hard-thresholding algorithms for low-rank matrix completion} \cite{jain2010guaranteed,vu2019local,vu2019accelerating}. 
In a recent work, Gratton and Tshimanga \cite{gratton2016second} presented a second-order expansion for the singular subspace decomposition and make use of the result to deduce the second-order sensitivity of the TSVD solution to least-squares problems. However, since their application focuses on the expansion of the truncated pseudo-inverse rather than the TSVD itself, no specific result in perturbation expansion of the TSVD is mentioned. 
\hlnew{In a different approach to analyzing the TSVD operator, Feppon and Lermusiaux }\cite{feppon2018geometric} \hlnew{studied the embedded geometry of the fixed-rank matrix manifold and characterized the projection onto it as a smooth ($C^\infty$) map.} Based on this geometric interpretation, the authors provided an explicit expression for the directional derivative of the TSVD of order $r$ at a certain matrix with \hlnew{rank greater than or equal to $r$.}\footnote{Despite the fact that Theorem~25 in \cite{feppon2018geometric} reads ``greater than $r$'', both the proof of the theorem and the direct communication with the authors (on September 17, 2020) suggest the result should also include the case of rank-$r$ matrices.} 
\hlnew{On the one hand, the result directly suggests the first-order perturbation expansion of the TSVD. On the other hand, the differential geometry-based approach, while offering a clear path for calculating the derivatives, does not offer a direct recipe for obtaining the error bound on the first-order approximation or the higher order terms in the expansion. At the time of writing this manuscript, we are not aware of any explicit expression of the second-order derivative of the TSVD.}

\hlnew{In this paper, we present a novel perturbation analysis of the truncated singular value decomposition. First, by utilizing the perturbation expansion for singular subspaces in} \cite{gratton2016second}, \hlnew{we derive the first-order perturbation expansion of the TSVD. Our result matches the result on the directional derivative of the TSVD in} \cite{feppon2018geometric}.
Furthermore, we extend our analysis to study the second-order perturbation expansion and show that when the matrix has exact rank $r$, the TSVD of order $r$ admits a simple expression for its second-order expansion. To the best of our knowledge, this is the first explicit result for the second-order perturbation expansion of the TSVD. \hlnew{Third, we establish an error bound on the first-order approximation of the TSVD about a rank-$r$ matrix. Our bound holds universally for any level (or magnitude) of the perturbation. Finally, we demonstrate how the proposed perturbation expansions and error bounds can be applied to study the mean squared error associated with the TSVD-based matrix denoising solution.}


\section{Notation and definitions}
\label{sec:def}

Throughout the paper, we use $\norm{\cdot}_F$ and $\norm{\cdot}_2$ to denote the Frobenius norm and the spectral norm of a matrix, respectively. Occasionally, $\norm{\cdot}_2$ is used on a vector to denote the Euclidean norm. 
Boldfaced symbols are reserved for vectors and matrices. In addition, the $s \times t$ all-zero matrix is denoted by $\bm 0_{s \times t}$ and the $s \times s$ identity matrix is denoted by $\bm I_s$. We also use $\bm e_i^s$ to denote the $i$-th vector in the natural basis of $\R^s$. When understood clearly from the context, the dimensions of vectors/matrices in the aforementioned notation may be omitted. As a slight abuse of notation, we define the big O notation for matrices as follows.
\begin{definition} \label{def:O}
Let $\bm \Delta$ be some matrix and $\bm F(\bm \Delta)$ be a matrix-valued function of $\bm \Delta$. Then, for any positive number $k$, $\bm F(\bm \Delta) = \bm \O(\norm{\bm \Delta}_F^k)$ if there exists some constant $0 \leq c < \infty$ such that 
\begin{align*}
    \lim_{\epsilon \to 0^+} \sup_{\norm{\bm \Delta}_F=\epsilon} \frac{\norm{\bm F(\bm \Delta)}_F}{\norm{\bm \Delta}_F^k} = c .
\end{align*}
\end{definition}
We emphasize the difference between the commonly used big O notation in the literature and the $\bm \O$ notation used in this manuscript. While the former requires $c$ to be strictly greater than $0$, our notation includes the case $c=0$ to imply both situations that $\bm F(\bm \Delta)$ approaches $\bm 0$ at a rate either equal or faster than $\norm{\bm \Delta}_F^k$. 
Similarly, when used for a vector, we replace the Frobenius norm by the Euclidean norm in Definition~\ref{def:O} to denote the corresponding quantity. 

In the rest of the paper, unless otherwise specified, the symbol $\bm X$ is used to denote an arbitrary matrix in $\R^{m \times n}$. Here, without loss of generality, we assume that $m \geq n$.
The SVD of $\bm X$ is written as $\bm X = \bm U \bm \Sigma \bm V^T$ where $\bm \Sigma$ is a $m \times n$ rectangular diagonal matrix with main diagonal entries are the singular values $\sigma_{1} \geq \sigma_{2} \geq \ldots \geq \sigma_{n} \geq 0$. For completeness, we denote the ``ghost'' singular values $\sigma_{n+1}=\ldots=\sigma_m=0$ in the case $m > n$. Additionally, $\bm U \in \R^{m \times m}$ and $\bm V \in \R^{n \times n}$ are orthogonal matrices such that $\bm U \bm U^T = \bm U^T \bm U = \bm I_m$ and $\bm V \bm V^T = \bm V^T \bm V = \bm I_n$. We note that the left and right singular vectors of $\bm X$ are the columns of $\bm U$ and $\bm V$, i.e., $\bm U = [\bm u_{1}, \bm u_{2}, \ldots, \bm u_{m}]$ and $\bm V = [\bm v_{1}, \bm v_{2}, \ldots, \bm v_{n}]$. 
Thus, $\bm X$ can also be rewritten as the sum of rank-$1$ matrices: $\bm X = \sum_{i=1}^n \sigma_{i} \bm u_{i} \bm v_{i}^T$.
Next, we define the singular subspace decomposition as follows. 
\begin{definition}
\label{def:SVD}
Given $1 \leq r < n$, the singular subspace decomposition of $\bm X \in \R^{m \times n}$ is given by:
\begin{align} \label{equ:X}
\bm X = \begin{bmatrix*}[r]
\bm U_{1} & \bm U_{2}
\end{bmatrix*} \begin{bmatrix*}
\bm \Sigma_{1} & \bm 0 \\
\bm 0 & \bm \Sigma_{2}
\end{bmatrix*} \begin{bmatrix*}[r]
\bm V_{1}^T \\ \bm V_{2}^T
\end{bmatrix*} = \bm U_{1} \bm \Sigma_{1} \bm V_{1}^T + \bm U_{2} \bm \Sigma_{2} \bm V_{2}^T ,
\end{align}
where
\begin{align*}
\bm \Sigma_{1} = \diag(\sigma_{1},\ldots,\sigma_{r}) \in \R^{r \times r}, \qquad \bm \Sigma_{2} = \begin{bmatrix*} \diag(\sigma_{r+1},\ldots,\sigma_{n}) \\ \bm 0 \end{bmatrix*}  \in \R^{(m-r) \times (n-r)},
\end{align*}
with the singular values in descending order, i.e., $\sigma_{1} \geq \sigma_{2} \geq \ldots \geq \sigma_{n} \geq 0$, and 
\begin{align*}
    &\bm U_{1} = \begin{bmatrix*}[r] \bm u_{1}& \ldots& \bm u_{r} \end{bmatrix*} \in \R^{m \times r}, \quad \bm U_{2} = \begin{bmatrix*}[r] \bm u_{r+1}&\ldots&\bm u_{m} \end{bmatrix*} \in \R^{m \times (m-r)}, \\
    &\bm V_{1} = \begin{bmatrix*}[r] \bm v_{1}& \ldots& \bm v_{r} \end{bmatrix*} \in \R^{n \times r}, \quad \bm V_{2} = \begin{bmatrix*}[r] \bm v_{r+1}&\ldots&\bm v_{n} \end{bmatrix*} \in \R^{n \times (n-r)} .
\end{align*}
\end{definition}
\noindent It is clear from Definition~\ref{def:SVD} that
\begin{align*}
    \bm U = \begin{bmatrix*}[r]
    \bm U_{1} & \bm U_{2}
    \end{bmatrix*}, \qquad \bm \Sigma = \begin{bmatrix*}
    \bm \Sigma_{1} & \bm 0 \\
    \bm 0 & \bm \Sigma_{2}
    \end{bmatrix*}, \qquad \bm V = \begin{bmatrix*}[r]
    \bm V_{1} & \bm V_{2}
    \end{bmatrix*} .
\end{align*}
Here the columns of $\bm U_1$ and $\bm U_2$ (or $\bm V_1$ and $\bm V_2$) provide the bases for the column-space (or row-space) of $\bm X$ and its orthogonal complement, respectively. 

\begin{definition} \label{def:Puv}
The orthonormal projectors onto the subspaces of $\bm X$ are defined as:
\begin{alignat*}{3}
    &\bm P_{\bm U_{1}} = \bm U_{1} \bm U_{1}^T = \sum_{i=1}^r \bm u_{i} \bm u_{i}^T , \qquad &&\bm P_{\bm U_{2}} = \bm U_{2} \bm U_{2}^T = \bm I_m - \bm P_{\bm U_{1}} = \sum_{i=r+1}^m \bm u_{i} \bm u_{i}^T , \\
    &\bm P_{\bm V_{1}} = \bm V_{1} \bm V_{1}^T = \sum_{i=1}^r \bm v_{i} \bm v_{i}^T , \qquad &&\bm P_{\bm V_{2}} = \bm V_{2} \bm V_{2}^T = \bm I_n - \bm P_{\bm V_{1}} = \sum_{i=r+1}^n \bm v_{i} \bm v_{i}^T . 
\end{alignat*}
\end{definition}
Generally,  matrices $\bm U_{1}$, $\bm U_{2}$, $\bm V_{1}$, and $\bm V_{2}$ are not unique. In particular, for simple non-zero singular values, the corresponding left and right singular vectors are unique up to a simultaneous sign change. For repeated and positive singular values, the corresponding left and right singular vectors are unique up to a simultaneous right multiplication with the same orthogonal matrix.  Finally, for zero singular values, the singular vectors can be any orthonormal bases of the left and right null spaces of $\bm X$. On the other hand, the singular subspaces spanned by the columns of $\bm U_{1}$, $\bm U_{2}$, $\bm V_{1}$, $\bm V_{2}$, and their corresponding projectors are unique provided that $\sigma_{r} > \sigma_{r+1}$ \cite{hansen1987truncatedsvd}. We are now in position to define the singular value truncation.
\begin{definition} \label{def:Pr}
The $r$-truncated singular value decomposition of $\bm X$ ($r$-TSVD) is defined as  
\begin{align} \label{equ:Eckart}
    \P_r(\bm X) = \sum_{i=1}^r \sigma_{i} \bm u_{i} \bm v_{i}^T = \bm U_{1} \bm \Sigma_{1} \bm V_{1}^T .
\end{align}
\end{definition}
By Eckart-Young theorem \cite{eckart1936approximation}, $\P_r(\bm X)$ is the best least squares approximation of $\bm X$ by a rank-$r$ matrix, with
respect to unitarily-invariant norms. Therefore, this operator is also known as the projection of $\bm X$ onto the non-convex set of rank-$r$ matrices.
$\P_r(\bm X)$ is unique if either $\sigma_{r} > \sigma_{r+1}$ or $\sigma_{r}=0$. 
In the special case when $\bm X$ has exact rank $r$, we have $\sigma_r > \sigma_{r+1} = \ldots = \sigma_n = 0$ and the projectors onto the subspaces of $\bm X$, namely, $\bm P_{\bm U_1},\bm P_{\bm U_2},\bm P_{\bm V_1},$ and $\bm P_{\bm V_2}$ are unique. However, the matrices $\bm U_2$ and $\bm V_2$ can take any orthonormal basis in $\bm R^{m-r}$ and $\bm R^{n-r}$, respectively, as their columns.  
Finally, for a rank-$r$ matrix, we define the pseudo inverse of $\bm X$ as $\bm X^\dagger = \bm U_1 \bm \Sigma_1^{-1} \bm V_1^T$. It is worth mentioning that $\norm{\bm X}_2 = \sigma_1$ while $\norm{\bm X^\dagger}_2=1/\sigma_r$ in this case.

\section{Preliminaries}
\label{sec:prel}

Two elemental bounds for singular values were given by Weyl \cite{weyl1912asymptotische} in 1912 and Mirsky \cite{mirsky1960symmetric} in 1960:
\begin{proposition}
\label{pro:s}
Let $\bm \Delta \in \R^{m \times n}$ be a perturbation of arbitrary magnitude. Denote $\tilde{\bm X} = \bm X + \bm \Delta$ with singular values $\tilde{\sigma}_{1} \geq \tilde{\sigma}_{2} \geq \ldots \geq \tilde{\sigma}_{n} \geq 0$. Then,
\begin{itemize}
    \item Weyl's inequality: $\abs{\tilde{\sigma}_{i} - \sigma_{i}} \leq \norm{\bm \Delta}_2$, for $i=1,\ldots,n$,
    \item Mirsky's inequality: $\sqrt{\sum_{i=1}^n (\tilde{\sigma}_{i} - \sigma_{i})^2} \leq \norm{\bm \Delta}_F$.
\end{itemize}
\end{proposition}
\noindent Proposition~\ref{pro:s} asserts that the changes in the singular values can \hlnew{be} bounded using only the norm of the perturbation.
By leveraging the specific values of the entries of the perturbation matrix, the behavior of singular values under perturbations can be described more precisely through perturbation expansions. In \cite{stewart1984second}, Stewart showed that if $\sigma_{n}$ is non-zero and distinct from other singular values of $\bm X$, then its corresponding perturbed singular value can be expressed by
\begin{align} \label{equ:ds1}
    \tilde{\sigma}_{n} = \sigma_{n} + \bm u_{n}^T \bm \Delta \bm v_{n} + \bm \O(\norm{\bm \Delta}^2) .
\end{align}
It is later known that the result in (\ref{equ:ds1}) also holds for any simple non-zero singular values \cite{Stewart90}. In another approach, Sun \cite{ji1988note} derived a second-order perturbation expansion for simple non-zero singular values.
For a simple zero singular value, Stewart \cite{stewart1984second} claimed that deriving a perturbation expansion is non-trivial and proposed a second-order approximation for $\tilde{\sigma}_n^2$ instead.
Most recently, a generalization of (\ref{equ:ds1}) to a set of singular values that is well separated from the rest is proved in \cite{stewart2006perturbation}.

While the singular values of a matrix are proven to be quite stable under perturbations, the singular vectors, especially those correspond to a cluster of singular values, are extremely sensitive. It is therefore natural to bound the perturbation error based on the subspace spanned by the singular vectors. Consider the singular subspace decomposition in Definition~\ref{def:SVD}. We define the singular gap as the smallest distance between a singular value in $\bm \Sigma_{1}$ and a singular value in $\bm \Sigma_{2}$.
When the spectral norm of the perturbation is smaller than this gap, Wedin's $\sin \Theta$ theorem \cite{wedin1972perturbation} provides an upper bound on the distances between the left and right singular subspaces and their corresponding perturbed counterparts in terms of the singular gap and the Frobenius norm of the perturbation. Furthermore, Stewart \cite{stewart1973error} showed that there exist explicit expressions of the perturbed subspaces in the bases of the unperturbed subspaces, which can be leveraged to obtain error bounds for certain characteristic subspaces associated with the SVD. Let us rephrase this result in the following proposition.
\begin{proposition} \label{pro:stewart} (Rephrased from Theorem~2.1 in \cite{gratton2016second}, which is based on Theorem~6.4 in \cite{stewart1973error})
In addition to the setting in Definition~\ref{def:SVD}, assume that $\sigma_r > \sigma_{r+1}$. For a perturbation $\bm \Delta \in \R^{m \times n}$, denote the singular subspace decomposition of $\tilde{\bm X} = \bm X + \bm \Delta$ by
\begin{align*}
    \tilde{\bm X} = \tilde{\bm U} \tilde{\bm \Sigma} \tilde{\bm V}^T = \begin{bmatrix*}[r] \tilde{\bm U}_1 & \tilde{\bm U}_2 \end{bmatrix*} \begin{bmatrix*} \tilde{\bm \Sigma}_1 & \bm 0 \\ \bm 0 & \tilde{\bm \Sigma}_2 \end{bmatrix*} \begin{bmatrix*}[r] \tilde{\bm V}_1^T \\ \tilde{\bm V}_2^T \end{bmatrix*} .
\end{align*}
Let us partition $\bm U^T \bm \Delta \bm V$ conformally with $\bm U$ and $\bm V$ in the form
\begin{align} \label{equ:E}
    \bm U^T \bm \Delta \bm V  = \begin{bmatrix*}[r] \bm U_1^T \bm \Delta \bm V_1 & \bm U_1^T \bm \Delta \bm V_2 \\ \bm U_2^T \bm \Delta \bm V_1 & \bm U_2^T \bm \Delta \bm V_2 \end{bmatrix*} = \begin{bmatrix*}[r] \bm E_{11} & \bm E_{12} \\ \bm E_{21} & \bm E_{22} \end{bmatrix*} = \bm E.
\end{align}
If
\begin{align} \label{cond:d}
    \norm{\bm \Delta}_2 < \frac{\sigma_r-\sigma_{r+1}}{2} ,
\end{align}
then there must exist unique matrices $\bm Q \in \R^{(m-r) \times r}$, $\bm P \in \R^{(n-r) \times r}$ whose norms are in the order of $\norm{\bm \Delta}_F$ such that
\begin{subequations}
\begin{align}
    &\bm Q (\bm \Sigma_1 + \bm E_{11}) + (\bm \Sigma_2 + \bm E_{22}) \bm P = - \bm E_{21} - \bm Q \bm E_{12} \bm P , \label{equ:QS} \\
    &(\bm \Sigma_1 + \bm E_{11}) \bm P^T + \bm Q^T (\bm \Sigma_2 + \bm E_{22}) = \bm E_{12} + \bm Q^T \bm E_{21} \bm P^T . \label{equ:PS}
\end{align} 
\label{equ:QPS}%
\end{subequations}
Moreover, using 
\begin{subequations}
\begin{align}
    \hat{\bm U}_1 &= (\bm U_1 - \bm U_2 \bm Q)(\bm I_r + \bm Q^T \bm Q)^{-1/2} , \label{equ:perturbed_U1} \\
    \hat{\bm U}_2 &= (\bm U_2 + \bm U_1 \bm Q^T)(\bm I_{m-r} + \bm Q \bm Q^T)^{-1/2} , \label{equ:perturbed_U2} \\
    \hat{\bm V}_1 &= (\bm V_1 + \bm V_2\bm P)(\bm I_r + \bm P^T\bm P)^{-1/2} , \label{equ:perturbed_V1} \\
    \hat{\bm V}_2 &= (\bm V_2 - \bm V_1 \bm P^T)(\bm I_{n-r} + \bm P \bm P^T)^{-1/2} , \label{equ:perturbed_V2}
\end{align}
\label{equ:perturbed}%
\end{subequations}
we can define semi-orthogonal matrices $\hat{\bm U}_1$, $\hat{\bm U}_2$, $\hat{\bm V}_1$, and $\hat{\bm V}_2$ satisfying $\hat{\bm U}_1^T \hat{\bm U}_2 = \bm 0$ and $\hat{\bm V}_1^T \hat{\bm V}_2 = \bm 0$, which provide bases to the same unique subspaces of $\tilde{\bm U}_1$, $\tilde{\bm U}_2$, $\tilde{\bm V}_1$, and $\tilde{\bm V}_2$, respectively, i.e., $\bm P_{\hat{\bm U}_1} = \bm P_{\tilde{\bm U}_1}$, $\bm P_{\hat{\bm U}_2} = \bm P_{\tilde{\bm U}_2}$, $\bm P_{\hat{\bm V}_1} = \bm P_{\tilde{\bm V}_1}$, and $\bm P_{\hat{\bm V}_2} = \bm P_{\tilde{\bm V}_2}$.
\end{proposition}
\noindent It is important to note that $\hat{\bm U}_1$, $\hat{\bm U}_2$, $\hat{\bm V}_1$, and $\hat{\bm V}_2$ may differ from $\tilde{\bm U}_1$, $\tilde{\bm U}_2$, $\tilde{\bm V}_1$, and $\tilde{\bm V}_2$, respectively. However, their corresponding subspaces are identical. This result will be useful later when replacing $\bm P_{\tilde{\bm U}_1}$ and $\bm P_{\tilde{\bm V}_1}$ in the following version of the $r$-TSVD $\P_r(\tilde{\bm X}) = \bm P_{\tilde{\bm U}_1} \tilde{\bm X} \bm P_{\tilde{\bm V}_1}$ with $\bm P_{\hat{\bm U}_1}$ and $\bm P_{\hat{\bm V}_1}$.
The substitution allows us to write an explicit expression of the $r$-TSVD using $\bm \Delta$ and terms that are in order of $\norm{\bm \Delta}_F$ such as $\bm Q$ and $\bm P$.
Equation (\ref{equ:QPS}) also enables the perturbation expansion of the SVD through the coefficient matrices $\bm Q$ and $\bm P$. 
In 1991, Li and Vaccaro \cite{li1991unified} considered a special case of rank-$r$ matrices ($\bm \Sigma_2 = \bm 0$) and introduced the first-order perturbation expansion for $\bm Q$ and $\bm P$ as a method to analyze the performance of subspace-based algorithms in array signal processing. Later on, Vaccaro \cite{vaccaro1994second} extended their approach to study the second-order perturbation expansion for the singular subspace decomposition. A more general result in this approach was proposed by Xu \cite{xu2002perturbation} in 2002, through relaxing the constraint $\bm \Sigma_2 = \bm 0$ to $\bm \Sigma_2^T \bm \Sigma_2 = \epsilon^2 \bm I$, for small $\epsilon \geq 0$. It was not until recently the second-order analysis with no restriction on $\bm \Sigma_2$ was provided by Gratton \cite{gratton2016second}. We summarize this result on second-order perturbation expansion for $\bm Q$ and $\bm P$ as follows.
\begin{proposition}
\label{pro:QP}
Given the setting in Proposition~\ref{pro:stewart}. Then
\begin{align} \label{equ:Q}
    \vect(\bm Q) = \bm \Phi_0^{-1} \bm \mu_1 + \bm \Phi_0^{-1} \bm \mu_2 - \bm \Phi_0^{-1} \bm \Phi_1 \bm \Phi_0^{-1} \bm \mu_1 + \bm \O(\norm{\bm \Delta}_F^3) ,
\end{align}
where 
\begin{align*}
    &\bm \Phi_0 = \bm \Sigma_1^2 \otimes \bm I_{m-r} - \bm I_r \otimes (\bm \Sigma_2 \bm \Sigma_2^T) , && \bm \Phi_1 = (\bm \Sigma_1 \bm E_{11}^T + \bm E_{11} \bm \Sigma_1) \otimes \bm I_{m-r} - \bm I_r \otimes (\bm \Sigma_2 \bm E_{22}^T + \bm E_{22} \bm \Sigma_2^T) ,\\
    &\bm \mu_1 = -\vect(\bm \Sigma_2 \bm E_{12}^T + \bm E_{21} \bm \Sigma_1), && \bm \mu_2 = -\vect(\bm E_{22} \bm E_{12}^T + \bm E_{21} \bm E_{11}^T) ,
\end{align*}
and
\begin{align} \label{equ:P}
    \vect(\bm P) = \bm \Psi^{-1} \bm \tau_1 + \bm \Psi_0^{-1} \bm \tau_2 - \bm \Psi_0^{-1} \bm \Psi_1 \bm \Psi_0^{-1} \bm \tau_1 + \bm \O(\norm{\bm \Delta}_F^3) ,
\end{align}
where
\begin{align*}
    &\bm \Psi_0 = \bm \Sigma_1^2 \otimes \bm I_{m-r} - \bm I_r \otimes (\bm \Sigma_2^T \bm \Sigma_2) , && \bm \Psi_1 = (\bm \Sigma_1 \bm E_{11} + \bm E_{11}^T \bm \Sigma_1) \otimes \bm I_{m-r} - \bm I_r \otimes (\bm \Sigma_2^T \bm E_{22} + \bm E_{22}^T \bm \Sigma_2) ,\\
    &\bm \tau_1 = -\vect(\bm \Sigma_2^T \bm E_{21} + \bm E_{12}^T \bm \Sigma_1), && \bm \tau_2 = -\vect(\bm E_{22}^T \bm E_{21} + \bm E_{12}^T \bm E_{11}) .
\end{align*}
\end{proposition}

\begin{corollary} \label{cor:QP}
Suppose in Proposition~\ref{pro:stewart}, $\bm X$ has rank $r$, i.e., $\bm \Sigma_2 = \bm 0$. Then
\begin{align*}
    &\bm Q = - \bm E_{21} \bm \Sigma_1^{-1} - \bm E_{22} \bm E_{12}^T \bm \Sigma_1^{-2} + \bm E_{21} \bm \Sigma_1^{-1} \bm E_{11} \bm \Sigma_1^{-1} + \bm \O(\norm{\bm \Delta}_F^3) , \\
    &\bm P = \bm E_{12}^T \bm \Sigma_1^{-1} + \bm E_{22}^T \bm E_{21} \bm \Sigma_1^{-2} - \bm E_{12}^T \bm \Sigma_1^{-1} \bm E_{11}^T \bm \Sigma_1^{-1} + \bm \O(\norm{\bm \Delta}_F^3) .
\end{align*}
\end{corollary}
Finally, we devote the rest of this section to discuss condition (\ref{cond:d}) in Proposition~\ref{pro:stewart}. As mentioned earlier, the singular subspaces corresponding to $\tilde{\bm U}_1$, $\tilde{\bm U}_2$, $\tilde{\bm V}_1$, and $\tilde{\bm V}_2$ are unique if and only if $\tilde{\sigma}_r > \tilde{\sigma}_{r+1}$. By Weyl's inequality (see Proposition~\ref{pro:s}), we have $\abs{\tilde{\sigma}_{r+1} - \sigma_{r+1}} \leq \norm{\bm \Delta}_2$. Since $\norm{\bm \Delta}_2  < (\sigma_r - \sigma_{r+1})/2$ and $\abs{\tilde{\sigma}_{r+1} - \sigma_{r+1}} \geq \tilde{\sigma}_{r+1} - \sigma_{r+1}$, one can further upper bound the $r+1$-th perturbed singular value by
\begin{align} \label{equ:w1}
     \tilde{\sigma}_{r+1} < \sigma_{r+1} + \frac{\sigma_r - \sigma_{r+1}}{2} = \frac{\sigma_r + \sigma_{r+1}}{2} .
\end{align}
Following a similar argument, $\abs{\tilde{\sigma}_{r} - \sigma_{r}} \leq \norm{\bm \Delta}_2$ leads to 
\begin{align} \label{equ:w2}
    \tilde{\sigma}_r > \sigma_r - \frac{\sigma_r - \sigma_{r+1}}{2} = \frac{\sigma_r + \sigma_{r+1}}{2} . 
\end{align}
It follows from (\ref{equ:w1}) and (\ref{equ:w2}) that the gap between $\tilde{\sigma}_r$ and $\tilde{\sigma}_{r+1}$ is strictly greater than $0$:
\begin{align} \label{equ:uniq}
    \tilde{\sigma}_{r+1} < \frac{\sigma_r + \sigma_{r+1}}{2} < \tilde{\sigma}_r .
\end{align}
As mentioned in \cite{gratton2016second}, condition (\ref{cond:d}) is more restrictive, but simpler, than the original condition specified in \cite{stewart1973error}. Based on the aforementioned preliminaries, we are ready to present our results.

\section{Perturbation expansions for the $r$-TSVD}
\label{sec:exp}

This section presents perturbation expansion results for the $r$-TSVD operator. In order to guarantee the uniqueness of the expansions, we assume \hlnew{throughout the section} that the $r$-th and $r+1$-th singular values are well-separated and the perturbation $\bm \Delta$ has small magnitude relative to $\bm X$. 

\hlnew{Let us begin with a non-trivial result on the first-order perturbation expansion of the $r$-TSVD. The result is consistent with Theorem~25 from} \cite{feppon2018geometric}, \hlnew{in which Feppon and Lermusiaux utilized differential geometry to derive a closed-form expression for the directional derivative of the $r$-TSVD. Using tools from perturbation analysis, we are able to obtain the same result on the first-order perturbation expansion of $\P_r$.} \hlnew{The additional benefit of the technique used here, as can be seen later, is that it can be leveraged to further derive the second-order perturbation expansion and the bound on the approximation error of the first-order expansion about a rank-$r$ matrix.}
\begin{theorem} \label{THEO:1ST}
Assume $\sigma_{r} > \sigma_{r+1}$. Then, for some perturbation $\bm \Delta \in \R^{m \times n}$ such that $\norm{\bm \Delta}_2 < \frac{\sigma_{r}-\sigma_{r+1}}{2}$, the first-order perturbation expansion of the $r$-TSVD about $\bm X$ is uniquely given by\footnote{\hlnew{We recall that throughout this manuscript we assume $m \geq n$.}}
\begin{align*}
    \P_r(\bm X + \bm \Delta) = \P_r(\bm X) + \bm \Delta - \bm P_{\bm U_{2}} \bm \Delta \bm P_{\bm V_{2}} + \sum_{i=1}^r &\sum_{j=r+1}^n \biggl( \frac{\sigma_{j}^2}{\sigma_{i}^2-\sigma_{j}^2} (\bm u_{i} \bm u_{i}^T \bm \Delta \bm v_{j} \bm v_{j}^T + \bm u_{j} \bm u_{j}^T \bm \Delta \bm v_{i} \bm v_{i}^T) \\
    &+ \frac{\sigma_{i} \sigma_{j}}{\sigma_{i}^2-\sigma_{j}^2} (\bm u_{i} \bm v_{i}^T \bm \Delta^T \bm u_{j} \bm v_{j}^T + \bm u_{j} \bm v_{j}^T \bm \Delta^T \bm u_{i} \bm v_{i}^T) \biggr) + \bm \O(\norm{\bm \Delta}_F^2) . \numberthis \label{equ:1st}
\end{align*}
\end{theorem}
\noindent The proof of Theorem~\ref{THEO:1ST} is based on perturbation expansions of the coefficient matrices $\bm Q$ and $\bm P$ in Proposition~\ref{pro:QP}. Interested readers are encouraged to find out the details in Appendix~\ref{apdx:1st}.
\hlnew{As mentioned earlier, the first-order term in} (\ref{equ:1st}) \hlnew{is equivalent to the directional derivative given by Theorem~25 in} \cite{feppon2018geometric}: \begin{align*}
    \bm \nabla_{\bm \Delta} \P_r(\bm X) = \bm P_{\bm U_2} \bm \Delta &\bm P_{\bm V_1} + \bm P_{\bm U_1} \bm \Delta \\
    &+ \sum_{i=1}^r \sum_{j=r+1}^{m} \frac{\sigma_j}{\sigma_i^2-\sigma_j^2} \biggl( \bigl( \sigma_i \bm u_j^T \bm \Delta \bm v_i + \sigma_j \bm u_i^T \bm \Delta \bm v_j \bigr) \bm u_j \bm v_i^T + \bigl( \sigma_j \bm u_j^T \bm \Delta \bm v_i + \sigma_i \bm u_i^T \bm \Delta \bm v_j \bigr) \bm u_i \bm v_j^T \biggr) . \numberthis \label{equ:Feppon}
\end{align*}
\hlnew{It is worthwhile to mention that we arrive at the first-order perturbation expansion in} Theorem~\ref{THEO:1ST} \hlnew{while working independently on the error bounds for TSVD} (see Section~\ref{sec:bound}). 

Note that the condition $\norm{\bm \Delta}_2 < (\sigma_{r}-\sigma_{r+1})/2$ guarantees a non-zero gap between the $r$-th and the $r+1$-th singular values of the perturbed matrix (see (\ref{equ:uniq})), and hence guarantees $\P_r(\bm X + \bm \Delta)$ \hlnew{on the LHS of} (\ref{equ:1st}) is unique.
\hlnew{At the same time}, each term on the RHS of (\ref{equ:1st}) is well-defined due to the uniqueness of singular subspaces associated with each group of singular values of $\bm X$. \hlnew{The term $\bm \Delta - \bm P_{\bm U_{2}} \bm \Delta \bm P_{\bm V_{2}}$ can be viewed as the projection of $\bm \Delta$ onto the tangent space of the manifold of rank-$r$ matrices }\cite{absil2012projection}. \hlnew{On the other hand, the double summation stems from the curvature of this manifold when $\bm X$ does not lie on it (with rank greater than $r$). To demonstrate the first-order expansion in} Theorem~\ref{THEO:1ST}, let us consider the following examples.

\begin{example} \label{eg:X1}
Consider the matrix $\bm X$ with its SVD as follows:
\begin{align} \label{equ:eg_X1}
    \bm X = \frac{1}{2} \begin{bmatrix*}[r] 4 & -4 & 7 \\ 0 & 0 & -9 \\ 4 & 8 & 1 \\ 8 & 4 & -1 \end{bmatrix*} = \left( \frac{1}{2} \begin{bmatrix*}[r] -1 & 1 & 1 & 1 \\ 1 & -1 & 1 & 1 \\ 1 & 1 & -1 & 1 \\ 1 & 1 & 1 & -1 \end{bmatrix*} \right) \cdot \begin{bmatrix*}[r] 6 & 0 & 0 \\ 0 & 6 & 0 \\ 0 & 0 & 3 \\ 0 & 0 & 0 \end{bmatrix*} \cdot \left( \frac{1}{3} \begin{bmatrix*}[r] 1 & 2 & 2 \\ 2 & 1 & -2 \\ -2 & 2 & -1 \end{bmatrix*} \right)^T .
\end{align}
In this example, note that $\sigma_1 = \sigma_2 > \sigma_3$. From Definition~\ref{def:Pr}, we have
\begin{align} \label{equ:eg_X2}
    \P_2(\bm X) = \left( \frac{1}{2} \begin{bmatrix*}[r] -1 & 1 \\ 1 & -1 \\ 1 & 1 \\ 1 & 1 \end{bmatrix*} \right) \cdot \begin{bmatrix*}[r] 6 & 0 \\ 0 & 6 \end{bmatrix*} \cdot \left( \frac{1}{3} \begin{bmatrix*}[r] 1 & 2 \\ 2 & 1 \\ -2 & 2 \end{bmatrix*} \right)^T = \begin{bmatrix*}[r] 1 & -1 & 4 \\ -1 & 1 & -4 \\ 3 & 3 & 0 \\ 3 & 3 & 0 \end{bmatrix*} .
\end{align}
In addition,
\begin{align} \label{equ:eg_X3}
    \bm P_{\bm U_2} = \frac{1}{2} \begin{bmatrix*}[r] 1 & 1 & 0 & 0 \\ 1 & 1 & 0 & 0 \\ 0 & 0 & 1 & -1 \\ 0 & 0 & -1 & 1 \end{bmatrix*}, \qquad \bm P_{\bm V_2} = \frac{1}{9} \begin{bmatrix*}[r] 4 & -4 & -2 \\ -4 & 4 & 2 \\ -2 & 2 & 1 \end{bmatrix*} .
\end{align}
For the perturbation
\begin{align*}
    \bm \Delta = \frac{3}{200} \begin{bmatrix*}[r] 3 & 3 & -9 \\ -3 & -9 & 3 \\ 7 & 5 & -5 \\ -1 & 7 & -7 \end{bmatrix*} , \text{ with } \norm{\bm \Delta}_F = 0.2985 < \delta = 1.5 , \numberthis \label{equ:eg_Delta}
\end{align*}
(\ref{equ:eg_X3}) leads to
\begin{align} \label{equ:eg_X4}
    \bm P_{\bm U_2} \bm \Delta \bm P_{\bm V_2} = \frac{3}{200} \begin{bmatrix*}[r] 2 & -2 & -1 \\ 2 & -2 & -1 \\ 2 & -2 & -1 \\ -2 & 2 & 1 \end{bmatrix*} .
\end{align}
Now the double summation in (\ref{equ:1st}) can be represented as
\begin{align*}
    \bm G(\bm \Delta) &= \frac{1}{3} (\bm u_1 \bm u_1^T \bm \Delta \bm v_3 \bm v_3^T + \bm u_3 \bm u_3^T \bm \Delta \bm v_1 \bm v_1^T) + \frac{2}{3} (\bm u_1 \bm v_1^T \bm \Delta^T \bm u_3 \bm v_3^T + \bm u_3 \bm v_3^T \bm \Delta^T \bm u_1 \bm v_1^T) \\
    &\qquad + \frac{1}{3} (\bm u_2 \bm u_2^T \bm \Delta \bm v_3 \bm v_3^T + \bm u_3 \bm u_3^T \bm \Delta \bm v_2 \bm v_2^T) + \frac{2}{3} (\bm u_2 \bm v_2^T \bm \Delta^T \bm u_3 \bm v_3^T + \bm u_3 \bm v_3^T \bm \Delta^T \bm u_2 \bm v_2^T) .
\end{align*}
While the singular vectors of $\bm X$ are not unique (due to $\sigma_1 = \sigma_2$), the singular subspaces of $\bm X$ are unique. Therefore, by representing $\bm G(\bm \Delta)$ as
\begin{align*}
    \bm G(\bm \Delta) &= \frac{1}{3} (\bm u_1 \bm u_1^T + \bm u_2 \bm u_2^T) \bm \Delta \bm v_3 \bm v_3^T + \frac{1}{3} \bm u_3 \bm u_3^T \bm \Delta (\bm v_1 \bm v_1^T + \bm v_2 \bm v_2^T) \\
    &\qquad + \frac{2}{3} (\bm u_1 \bm v_1^T + \bm u_2 \bm v_2^T) \bm \Delta^T \bm u_3 \bm v_3^T + \frac{2}{3} \bm u_3 \bm v_3^T \bm \Delta^T (\bm u_1 \bm v_1^T + \bm u_2 \bm v_2^T) , \numberthis \label{equ:egG}
\end{align*}
we observe that $\bm G(\bm \Delta)$ is well-defined since $\bm u_1 \bm u_1^T + \bm u_2 \bm u_2^T$, $\bm u_3 \bm u_3^T$, $\bm v_1 \bm v_1^T + \bm v_2 \bm v_2^T$, $\bm v_3 \bm v_3^T$, $\bm u_1 \bm v_1^T + \bm u_2 \bm v_2^T$, and $\bm u_3 \bm v_3^T$ are all unique quantities, namely,
\begin{align*}
    &\bm u_1 \bm u_1^T + \bm u_2 \bm u_2^T = \bm P_{\bm U_1} = \frac{1}{2} \begin{bmatrix*}[r] 1 & -1 & 0 & 0 \\ -1 & 1 & 0 & 0 \\ 0 & 0 & 1 & 1 \\ 0 & 0 & 1 & 1 \end{bmatrix*}, && \bm u_3 \bm u_3^T = \frac{1}{4} \begin{bmatrix*}[r] 1 & 1 & -1 & 1 \\ 1 & 1 & -1 & 1 \\ -1 & -1 & 1 & -1 \\ 1 & 1 & -1 & 1 \end{bmatrix*}, \\
    &\bm v_1 \bm v_1^T + \bm v_2 \bm v_2^T = \bm P_{\bm V_1} = \frac{1}{9} \begin{bmatrix*}[r] 5 & 4 & 2 \\ 4 & 5 & -2 \\ 2 & -2 & 8 \end{bmatrix*} , && \bm v_3 \bm v_3^T = \bm P_{\bm V_2} = \frac{1}{9} \begin{bmatrix*}[r] 4 & -4 & -2 \\ -4 & 4 & 2 \\ -2 & 2 & 1 \end{bmatrix*}, \\
    &\bm u_1 \bm v_1^T + \bm u_2 \bm v_2^T = \frac{1}{18} \begin{bmatrix*}[r] 3 & -3 & 12 \\ -3 & 3 & -12 \\ 9 & 9 & 0 \\ 9 & 9 & 0 \end{bmatrix*}, && \bm u_3 \bm v_3^T = \frac{1}{6} \begin{bmatrix*}[r] 2 & -2 & -1 \\ 2 & -2 & -1 \\ -2 & 2 & 1 \\ 2 & -2 & -1 \end{bmatrix*} . \numberthis \label{equ:uv_eg}
\end{align*}
Substituting the values of the $6$ aforementioned terms in (\ref{equ:uv_eg}) and the value of $\bm \Delta$ in (\ref{equ:eg_Delta}) back into (\ref{equ:egG}), we obtain 
\begin{align} \label{equ:eg_X5}
    \bm G(\bm \Delta) = \frac{1}{200} \begin{bmatrix*}[r] -6 & 3 & 0 \\ 2 & -5 & -4 \\ -2 & 5 & 4 \\ -6 & 3 & 0 \end{bmatrix*} .
\end{align}
The substitution of (\ref{equ:eg_X2}), (\ref{equ:eg_Delta}), (\ref{equ:eg_X4}), and (\ref{equ:eg_X5}) into (\ref{equ:1st}) yields
\begin{align} \label{equ:eg_X6}
    \P_2(\bm X + \bm \Delta) = \begin{bmatrix*}[r] 0.9850 & -0.9100 & 3.8800 \\ -1.0650 & 0.8700 & -3.9600 \\ 3.0600 & 3.1300 & -0.0400 \\ 2.9850 & 3.0900 & -0.1200 \end{bmatrix*} + \bm \O(\norm{\bm \Delta}_F^2) .
\end{align}
On the other hand, running a simple numerical evaluation by Definition~\ref{def:Pr}, we can compute $\P_2(\bm X + \bm \Delta)$ and obtain
\begin{align*}
    \P_2(\bm X + \bm \Delta) = \begin{bmatrix*}[r] 0.9840 & -0.9088 & 3.8792 \\ -1.0632 & 0.8689 & -3.9615 \\ 3.0650 & 3.1284 & -0.0403 \\ 2.9870 & 3.0890 & -0.1213 \end{bmatrix*} .
\end{align*}
The approximation error of the first-order perturbation expansion has magnitude of $0.0043$, which is much smaller than the approximation error of the zero-order expansion, i.e., $\norm{\P_2(\bm X + \bm \Delta) - \P_2(\bm X)}_F = 0.3016$.
\end{example}


\begin{example}
Let us consider a counter-example in which the condition $\norm{\bm \Delta}_2 < (\sigma_{r}-\sigma_{r+1})/2$ is not satisfied. In particular, by setting
\begin{align*}
    \bm X = \begin{bmatrix*}[r] 2 & 0 & 0 \\ 0 & 2 & 0 \\ 0 & 0 & 1 \\ 0 & 0 & 0 \end{bmatrix*} , \qquad \bm \Delta = \begin{bmatrix*} 0.1 & 0 & 0 \\ 0 & -0.5 & 0 \\ 0 & 0 & 0.5 \\ 0 & 0 & 0 \end{bmatrix*} ,
\end{align*}
following similar calculation in Example~\ref{eg:X1} would yield
\begin{align*}
    \P_2(\bm X + \bm \Delta) = \begin{bmatrix*} 2.1 & 0 & 0 \\ 0 & 1.5 & 0 \\ 0 & 0 & 0 \\ 0 & 0 & 0 \end{bmatrix*} + \bm \O(\norm{\bm \Delta}_F^2) .
\end{align*}
On the other hand, the $2$-TSVD of $\bm X + \bm \Delta$ can either be
\begin{align*}
    \P_2(\bm X + \bm \Delta) = \begin{bmatrix*} 2.1 & 0 & 0 \\ 0 & 1.5 & 0 \\ 0 & 0 & 0 \\ 0 & 0 & 0 \end{bmatrix*} \qquad \text{or} \qquad \P_2(\bm X + \bm \Delta) = \begin{bmatrix*} 2.1 & 0 & 0 \\ 0 & 0 & 0 \\ 0 & 0 & 1.5 \\ 0 & 0 & 0 \end{bmatrix*} .
\end{align*}
It can be seen that our first-order approximation is no longer accurate if the later truncation is considered.
\end{example}


One immediate consequence of Theorem~\ref{THEO:1ST} is when the matrix has exact rank $r$, the double summation on the RHS of (\ref{equ:1st}) vanishes since $\sigma_{j} = 0$ for all $j>r$. Thus, we obtain a simple expression for the first-order expansion of $\P_r(\cdot)$ about a rank-$r$ matrix.

\begin{corollary} \label{COR:1ST}
Let $\bm X \in \R^{m \times n}$ be a rank-$r$ matrix. Then, for some perturbation $\bm \Delta \in \R^{m \times n}$ such that $\norm{\bm \Delta}_2 < \sigma_r/2$, the first-order perturbation expansion of the $r$-TSVD about $\bm X$ is uniquely given by
\begin{align} \label{equ:1str}
    \P_r(\bm X+ \bm \Delta) = \bm X + \bm \Delta - \bm P_{\bm U_2} \bm \Delta \bm P_{\bm V_2} + \bm \O(\norm{\bm \Delta}_F^2) . 
\end{align}
\end{corollary}
\noindent We observe that while the first-order term depends on the perturbation $\bm \Delta$ and the two projections $\bm P_{\bm U_2}$ and $\bm P_{\bm V_2}$, it is independent of the singular values of $\bm X$. Motivated by the simple result in Corollary~\ref{COR:1ST}, we further study the second-order perturbation expansion of the $r$-TSVD about a rank-$r$ matrix in the following theorem. 
\begin{theorem}
\label{THEO:2ND}
Let $\bm X \in \R^{m \times n}$ be a rank-$r$ matrix. Then, for some perturbation $\bm \Delta \in \R^{m \times n}$ such that $\norm{\bm \Delta}_2 < \sigma_r/2$, the second-order perturbation expansion of the $r$-TSVD about $\bm X$ is uniquely given by
\begin{align} \label{equ:rankop2}
    \P_r (\bm X + \bm \Delta) = \bm X + \bm \Delta - \bm P_{\bm U_2} \bm \Delta \bm P_{\bm V_2} + \bm X^\dagger \bm \Delta^T \bm P_{\bm U_2} \bm \Delta \bm P_{\bm V_2} + \bm P_{\bm U_2} \bm \Delta \bm P_{\bm V_2} \bm \Delta^T \bm X^\dagger + \bm P_{\bm U_2} \bm \Delta {(\bm X^\dagger)}^T \bm \Delta \bm P_{\bm V_2} + \bm \O(\norm{\bm \Delta}_F^3) .
\end{align}
\end{theorem}
\noindent The proof of Theorem~\ref{THEO:2ND} is given in Appendix~\ref{apdx:2nd}. \hlnew{The theorem states} that $\P_r(\bm X + \bm \Delta)$ admits a simple second-order approximation that only depends on $\bm P_{\bm U_2}$, $\bm P_{\bm V_2}$, and $\bm X^\dagger$ in addition to $\bm X$ and $\bm \Delta$ themselves. Notice the dependence of the three second-order terms on the RHS of (\ref{equ:rankop2}) on the pseudo inverse of $\bm X$ \hlnew{indicates the first-order approximation is sensitive to the least singular value of $\bm X$. In the next section, we shall prove that the error bound for the first-order approximation of $\P_r(\bm X + \bm \Delta)$ depends linearly on $1/\sigma_r$.}

\begin{remark} \label{rem:2nd0}
\hlnew{The differentiability of $\P_r$ at a rank-$r$ matrix, as shown in Corollary~2 and Theorem~2, matches with the well-known result in differential geometry that a projection onto the base of the normal bundle of any smooth manifold is a smooth map on the tubular neighborhood} \cite{lee2013smooth}. \hlnew{In particular, $\P_r$ is a classic smooth ($C^\infty$) map in a small open neighborhood containing the manifold of rank-$r$ matrices.}
\end{remark}

\begin{remark} \label{rem:2nd}
\hlnew{It is known that} the $r$-TSVD is differentiable at any point (matrix) with a non-zero gap between the $r$-th and $r+1$-th singular values \hlnew{and hence, admits a first-order perturbation expansion about such point}. While our result in Theorem~\ref{THEO:2ND} only considers a special case of rank-$r$ matrices, we suspect there exists a second-order perturbation expansion of the $r$-TSVD about a matrix $\bm X$ with rank greater than $r$. However, given the complexity of the first-order expansion, it certainly requires more elaborate work. We leave this as a future research direction. 
\end{remark}

\section{Error bounds for the $r$-TSVD}
\label{sec:bound}

This section introduces upper bounds on the difference between the $r$-TSVD and its first-order approximation. While in Section~\ref{sec:exp} the perturbation expansions are derived under the assumption that $\norm{\bm \Delta}_2 < (\sigma_r - \sigma_{r+1})/2$, the error bounds in this section do not require this constraint and indeed they hold for $\bm \Delta$ with arbitrary magnitude. 
It is important to note that, without the constraint on the level of the perturbation, $\P_r (\bm X + \bm \Delta)$ may not be unique since there is no guarantee that $\tilde{\sigma}_r > \tilde{\sigma}_{r+1}$. The value of $\P_r (\bm X + \bm \Delta)$ in case $\tilde{\sigma}_r = \tilde{\sigma}_{r+1}$ depends on the choice of the singular subspace decomposition of $\tilde{\bm X} = \bm X + \bm \Delta$ (see Definition~\ref{def:SVD}). Nevertheless, we shall provide error bounds that hold independent of the choice of decomposition. 

Let us consider the first-order expansion in (\ref{equ:1str}). One trivial bound on the approximation error can be derived as follows (see details in Appendix~\ref{apdx:bound0}):
\begin{lemma} \label{lem:bound0}
Let $\bm X \in \R^{m \times n}$ be a rank-$r$ matrix. For any $\bm \Delta \in \R^{m \times n}$  and any valid choice of subspace decomposition of $\bm X + \bm \Delta$, we have
\begin{align*}
    \norm{\P_r (\bm X + \bm \Delta) - (\bm X + \bm \Delta - \bm P_{\bm U_2} \bm \Delta \bm P_{\bm V_2})}_F \leq \norm{\bm X}_F + 2\norm{\bm \Delta}_F .
\end{align*}
\end{lemma}
\noindent Lemma~\ref{lem:bound0} suggests that for large ${\bm \Delta}$, the approximation error grows at most linearly in the norm of $\bm \Delta$. However, for small $\bm \Delta$, the aforementioned bound is not tight since Corollary~\ref{COR:1ST} implies the error should be in the order of $\norm{\bm \Delta}_F^2$. In order to tighten the bound for the small perturbation, we need to develop a different approach that is more meticulous about intermediate inequalities. We state our main result regarding the global error bound on the first-order approximation of the $r$-TSVD as follows. 
\begin{theorem}
\label{THEO:RANKOP}
Let $\bm X \in \R^{m \times n}$ be a rank-$r$ matrix. Then, for any $\bm \Delta \in \R^{m \times n}$ and any valid choice of subspace decomposition of $\bm X + \bm \Delta$, the first-order Taylor expansion of the $r$-TSVD about $\bm X$ is given by
\begin{align} \label{equ:rankop}
\P_r (\bm X + \bm \Delta) = \bm X + \bm \Delta - \bm P_{\bm U_2} \bm \Delta \bm P_{\bm V_2} + \bm R_{\bm X}(\bm \Delta) ,
\end{align}
where \hlnew{there exists a universal constant $1 + 1/\sqrt{2} \leq c \leq 4(1+\sqrt{2})$ such that}
\begin{align} \label{equ:boundR_1}
    \norm{\bm R_{\bm X}(\bm \Delta)}_F \leq  \frac{c}{\sigma_r} \norm{\bm \Delta}_F^2 .
\end{align}
\hlnew{Furthermore, the following inequality holds}
\begin{align} \label{equ:boundR_2}
    \norm{\bm R_{\bm X}(\bm \Delta)}_F \leq 2(1+\sqrt{2}) \norm{\bm \Delta}_F \min \biggl\{ \frac{2}{\sigma_r} \norm{\bm \Delta}_F, 1 \biggr\} .
\end{align}
\end{theorem}
\noindent The proof of Theorem~\ref{THEO:RANKOP} is given in Appendix~\ref{apdx:rankop}. 
It is noticeable that the first three terms on the RHS of (\ref{equ:rankop}) are uniquely given by the \hlnew{rank-$r$} singular subspace decomposition of $\bm X$. On the contrary, the LHS may not be unique (e.g., when $\tilde{\sigma}_r = \tilde{\sigma}_{r+1}$) and hence, so does the residual $\bm R_{\bm X}(\bm \Delta)$. However, it is interesting to note that the theorem makes no assumption on the norm of $\bm \Delta$, as well as the choice of the r-TSVD of $\bm X + \bm \Delta$.
The bound on the residual (or the remainder) in Theorem~\ref{THEO:RANKOP} is similar to the Lagrange error bound in univariate first-order Taylor series. \hlnew{It not only asserts that the approximation error can grow no faster than a quadratic rate but also determines the constant attached to $\norm{\bm \Delta}_F^2$.} 
Furthermore, the bound depends only on the $\sigma_r$ and $\norm{\bm \Delta}_F$, as one may expect from the second-order perturbation expansion of the $r$-TSVD in Theorem~\ref{THEO:2ND}. 
\begin{remark} \label{rem:c}
We conjecture but are unable to prove that the lower bound on $c$ is tight, i.e., $c = 1 + 1/\sqrt{2}$. Partial result in this direction regarding $\bm \Delta$ of certain structure is also given in the proof of Theorem~\ref{THEO:RANKOP}. \hlnew{In our numerical experiment, we ran multiple optimization procedures to maximize the quantity $\sigma_r \norm{\bm R_{\bm X}(\bm \Delta)}_F / \norm{\bm \Delta}_F^2$ with respect to $\bm \Delta$ and obtained the same constant $1 + 1/\sqrt{2}$.}
\end{remark}


\noindent The bound in (\ref{equ:boundR_2}) suggests an interesting behavior of the residual $\bm R_{\bm X}(\bm \Delta)$. When the perturbation is small, the error depends quadratically on the magnitude of the perturbation and inversely proportional to the least singular value of $\bm X$. In particular, as $\sigma_r$ approaches $0$, the first-order approximation becomes less accurate. On the contrary, for large $\bm \Delta$, the upper bound is linear in the norm of $\bm \Delta$ and independent of $\sigma_r$. Compared to the bound in Lemma~\ref{lem:bound0}, we observe that the dependence on $\bm X$ is eliminated. Asymptotically as $\norm{\bm \Delta}_F$ approaches $\infty$, the simple bound in the lemma becomes tighter than the bound in (\ref{equ:boundR_2}). 

We conclude this section by describing the behavior of the residual term for small perturbations. 
While it is challenging to establish a tight bound on $\norm{\bm R_{\bm X}(\bm \Delta)}_F$ (as a function of $\norm{\bm \Delta}_F$) for large $\bm \Delta$, it is possible to project the first-order approximation error for small perturbation based on the knowledge of the second-order perturbation expansion of the $r$-TSVD (see Theorem~\ref{THEO:2ND}). 
We provide the result in the following theorem, with the proof given in Appendix~\ref{apdx:RSUP}.
\begin{theorem} \label{THEO:RSUP}
Asymptotically as $\norm{\bm \Delta}_F$ approaches $0$, the norm of the residual term in Theorem~\ref{THEO:RANKOP} can be upper-bounded tightly by
\begin{align*}
    \lim_{\epsilon \to 0^+} \sup_{\norm{\bm \Delta}_F = \epsilon} \frac{\norm{\bm R_{\bm X}(\bm \Delta)}_F}{\norm{\bm \Delta}_F^2} = \frac{1}{\sigma_r \sqrt{3}} .    
\end{align*}
\end{theorem}

\begin{remark}
While Theorem~\ref{THEO:1ST} provides the first-order perturbation expansion of the $r$-TSVD about an arbitrary matrix $\bm X$ with $\sigma_r > \sigma_{r+1} \geq 0$, extending Theorems~\ref{THEO:RANKOP} and \ref{THEO:RSUP} to that case remains to be one of our future research directions due to the difficulty of bounding the double summation in (\ref{equ:1st}). 
\end{remark}

\section{\hlnew{An application to performance analysis in matrix denoising}}
This section presents an application of our result to the performance analysis of the TSVD for matrix denoising.
In many applications such as image denoising \cite{guo2015efficient}, multi-input multi-output (MIMO) channel estimation \cite{lindskog1999reduced}, collaborative filtering \cite{koren2009matrix}, low-rank procedures are often motivated by the following statistical model:
\begin{align*}
    \tilde{\bm X} = \bm X + \bm \Delta ,
\end{align*}
where $\bm X \in \R^{m \times n}$ is the unknown matrix with rank $r \leq \min (m,n)$ and $\bm \Delta$ is a random matrix whose entries are $i.i.d.$ normally distributed with zero mean and $\sigma^2$-variance, i.e., $\Delta_{ij} \sim \N(0,\sigma^2)$ for $i=1,\ldots,m$ and $j=1,\ldots,n$. To denoise the data, the TSVD is applied to the noisy matrix $\tilde{\bm X}$ to obtain the following estimator:
\begin{align*}
    \hat{\bm X} = \P_r (\tilde{\bm X}) .
\end{align*}
We would like to assess the mean squared error (MSE) of this estimator using our perturbation analysis of the TSVD.
As a baseline for our analysis, we consider the MSE of the noisy matrix $\tilde{\bm X}$:
\begin{align*}
    \E\bigl[\norm{\tilde{\bm X} - \bm X}_F^2\bigr] = \E\bigl[ \norm{\bm \Delta}_F^2 \bigr] = \sum_{i=1}^m \sum_{j=1}^n \E \bigl[ \Delta_{ij}^2 \bigr] = \sigma^2 mn . \numberthis \label{equ:M_tilde}
\end{align*}
Next, we study the MSE of the estimator $\hat{\bm X}$, i.e., $\E [ \norm{\hat{\bm X} - \bm X}_F^2 ]$. To the best of our knowledge, there exists no closed-form expression of this quantity due to the non-linearity of the truncated singular value operator. In the following, we provide the first-order approximation, the second-order approximation, and the upper bound for $\E [ \norm{\hat{\bm X} - \bm X}_F^2 ]$ based on the results presented in this paper. 

\begin{enumerate}[leftmargin=*]
\item \textbf{The first-order approximation:}

Let $\hat{\bm X}_1 = \bm X + \bm \Delta - \bm P_{\bm U_2} \bm \Delta \bm P_{\bm V_2}$ be the first-order approximation of $\hat{\bm X}$. 
We have
\begin{align*}
    \E \Bigl[ \norm{\hat{\bm X}_1 - \bm X}_F^2 \Bigr] &= \E \Bigl[ \norm{\bm \Delta - \bm P_{\bm U_2} \bm \Delta \bm P_{\bm V_2}}_F^2 \Bigr] \\
    &= \E \Bigl[ \norm{(\bm I_{mn} - \bm P_{\bm V_2} \otimes \bm P_{\bm U_2}) \vect (\bm \Delta)}_2^2 \Bigr] && \text{(by Lemma~\ref{lem:kron}-2)} \\
    &= \E \Bigl[ \bigl(\vect(\bm \Delta) \bigr)^T (\bm I_{mn} - \bm P_{\bm V_2} \otimes \bm P_{\bm U_2})^T (\bm I_{mn} - \bm P_{\bm V_2} \otimes \bm P_{\bm U_2}) \vect(\bm \Delta) \Bigr] . \numberthis \label{equ:M1a}
\end{align*}
Using the fact that $\bm P_{\bm U_2}$ and $\bm P_{\bm V_2}$ are projection matrices,  $\bm I_{mn} - \bm P_{\bm V_2} \otimes \bm P_{\bm U_2}$ is also a projection matrix, and hence, $(\bm I_{mn} - \bm P_{\bm V_2} \otimes \bm P_{\bm U_2})^T (\bm I_{mn} - \bm P_{\bm V_2} \otimes \bm P_{\bm U_2}) = \bm I_{mn} - \bm P_{\bm V_2} \otimes \bm P_{\bm U_2}$. Thus, (\ref{equ:M1a}) can be further simplified as
\begin{align*}
    \E \Bigl[ \norm{\hat{\bm X}_1 - \bm X}_F^2 \Bigr] &= \E \Bigl[ \bigl(\vect(\bm \Delta) \bigr)^T \bigr) (\bm I_{mn} - \bm P_{\bm V_2} \otimes \bm P_{\bm U_2}) \vect(\bm \Delta) \Bigr] \\
    &= \E \Bigl[ \tr \bigl( (\bm I_{mn} - \bm P_{\bm V_2} \otimes \bm P_{\bm U_2}) \vect(\bm \Delta) \bigl(\vect(\bm \Delta) \bigr)^T \bigr) \Bigr] && \text{(by the cyclic property of the trace)} \\
    &= \tr \Bigl( (\bm I_{mn} - \bm P_{\bm V_2} \otimes \bm P_{\bm U_2}) \E\bigl[ \vect(\bm \Delta) \bigl(\vect(\bm \Delta) \bigr)^T \bigr] \Bigr) .
\end{align*}
Since $\Delta_{ij} \overset{i.i.d.}{\sim} \N(0,\sigma^2)$, $\E\bigl[ \vect(\bm \Delta) \bigl(\vect(\bm \Delta) \bigr)^T \bigr] = \sigma^2 \bm I_{mn}$. Thus,
\begin{align*}
    \E \Bigl[ \norm{\hat{\bm X}_1 - \bm X}_F^2 \Bigr] &= \sigma^2 \tr(\bm I_{mn} - \bm P_{\bm V_2} \otimes \bm P_{\bm U_2}) \\
    &= \sigma^2 \tr(\bm I_{mn}) - \tr(\bm P_{\bm V_2}) \tr(\bm P_{\bm U_2}) \\
    &= \sigma^2 r (m+n-r) , \numberthis \label{equ:M1}
\end{align*}
where the second equality uses Lemma~\ref{lem:kron}-4 and the third equality stems from the fact that $\tr(\bm P_{\bm U_2}) = \tr(\bm U_2 \bm U_2^T) = \tr(\bm U_2^T \bm U_2) = \tr(\bm I_{m-r}) = m-r$ (and similarly $\tr(\bm P_{\bm V_2}) = n-r$).

\item \textbf{The second-order approximation:}

Let $\hat{\bm X}_2 = \bm X + \bm \Delta - \bm P_{\bm U_2} \bm \Delta \bm P_{\bm V_2} + \bm X^\dagger \bm \Delta^T \bm P_{\bm U_2} \bm \Delta \bm P_{\bm V_2} + \bm P_{\bm U_2} \bm \Delta \bm P_{\bm V_2} \bm \Delta^T \bm X^\dagger + \bm P_{\bm U_2} \bm \Delta {(\bm X^\dagger)}^T \bm \Delta \bm P_{\bm V_2}$ be the second-order approximation of $\hat{\bm X}$.
We have
\begin{align*}
    \E \Bigl[ &\norm{\hat{\bm X}_2 - \bm X}_F^2 \Bigr] = \E \Bigl[ \norm{\bm \Delta - \bm P_{\bm U_2} \bm \Delta \bm P_{\bm V_2} + \bm X^\dagger \bm \Delta^T \bm P_{\bm U_2} \bm \Delta \bm P_{\bm V_2} + \bm P_{\bm U_2} \bm \Delta \bm P_{\bm V_2} \bm \Delta^T \bm X^\dagger + \bm P_{\bm U_2} \bm \Delta {(\bm X^\dagger)}^T \bm \Delta \bm P_{\bm V_2}}_F^2 \Bigr] \\
    &= \E \bigl[ \norm{\bm \Delta - \bm P_{\bm U_2} \bm \Delta \bm P_{\bm V_2}}_F^2 \bigr] + \E \bigl[ \norm{\bm X^\dagger \bm \Delta^T \bm P_{\bm U_2} \bm \Delta \bm P_{\bm V_2} + \bm P_{\bm U_2} \bm \Delta \bm P_{\bm V_2} \bm \Delta^T \bm X^\dagger + \bm P_{\bm U_2} \bm \Delta {(\bm X^\dagger)}^T \bm \Delta \bm P_{\bm V_2}}_F^2 \bigr] \\
    &\qquad \qquad + \E \Bigl[ 2\tr \bigl( (\bm \Delta - \bm P_{\bm U_2} \bm \Delta \bm P_{\bm V_2})^T (\bm X^\dagger \bm \Delta^T \bm P_{\bm U_2} \bm \Delta \bm P_{\bm V_2} + \bm P_{\bm U_2} \bm \Delta \bm P_{\bm V_2} \bm \Delta^T \bm X^\dagger + \bm P_{\bm U_2} \bm \Delta {(\bm X^\dagger)}^T \bm \Delta \bm P_{\bm V_2}) \bigr) \Bigr]. \numberthis \label{equ:M2nd}
\end{align*}
Since $\Delta_{ij} \overset{i.i.d.}{\sim} \N(0,\sigma^2)$, the expected value of the third-order term on the RHS of (\ref{equ:M2nd}) is zero, i.e., 
\begin{align*}
    \E \Bigl[ 2\tr \bigl( (\bm \Delta - \bm P_{\bm U_2} \bm \Delta \bm P_{\bm V_2})^T (\bm X^\dagger \bm \Delta^T \bm P_{\bm U_2} \bm \Delta \bm P_{\bm V_2} + \bm P_{\bm U_2} \bm \Delta \bm P_{\bm V_2} \bm \Delta^T \bm X^\dagger + \bm P_{\bm U_2} \bm \Delta {(\bm X^\dagger)}^T \bm \Delta \bm P_{\bm V_2}) \bigr) \Bigr] = 0 .
\end{align*}
Therefore,
\begin{align*}
    \E \Bigl[ \norm{\hat{\bm X}_2 - \bm X}_F^2 \Bigr] = \E \bigl[ \norm{\bm \Delta - \bm P_{\bm U_2} \bm \Delta \bm P_{\bm V_2}}_F^2 \bigr] + \E \bigl[ \norm{\bm X^\dagger \bm \Delta^T \bm P_{\bm U_2} \bm \Delta \bm P_{\bm V_2} + \bm P_{\bm U_2} \bm \Delta \bm P_{\bm V_2} \bm \Delta^T \bm X^\dagger + \bm P_{\bm U_2} \bm \Delta {(\bm X^\dagger)}^T \bm \Delta \bm P_{\bm V_2}}_F^2 \bigr] .
\end{align*}
Since the first term on the RHS is given by (\ref{equ:M1}), we proceed with the calculation of the second term on the RHS, i.e., $\E \bigl[ \norm{\bm X^\dagger \bm \Delta^T \bm P_{\bm U_2} \bm \Delta \bm P_{\bm V_2} + \bm P_{\bm U_2} \bm \Delta \bm P_{\bm V_2} \bm \Delta^T \bm X^\dagger + \bm P_{\bm U_2} \bm \Delta {(\bm X^\dagger)}^T \bm \Delta \bm P_{\bm V_2}}_F^2 \bigr]$. Since $\bm X^\dagger = \bm P_{\bm U_1} \bm X^\dagger \bm P_{\bm V_1}$, the three terms inside the norm are orthogonal to each other, i.e., their inner products are zero. Hence,
\begin{align*}
    &\norm{\bm X^\dagger \bm \Delta^T \bm P_{\bm U_2} \bm \Delta \bm P_{\bm V_2} + \bm P_{\bm U_2} \bm \Delta \bm P_{\bm V_2} \bm \Delta^T \bm X^\dagger + \bm P_{\bm U_2} \bm \Delta {(\bm X^\dagger)}^T \bm \Delta \bm P_{\bm V_2}}_F^2 \\
    &\qquad \qquad \qquad \qquad = \norm{\bm X^\dagger \bm \Delta^T \bm P_{\bm U_2} \bm \Delta \bm P_{\bm V_2}}_F^2 + \norm{\bm P_{\bm U_2} \bm \Delta \bm P_{\bm V_2} \bm \Delta^T \bm X^\dagger}_F^2 + \norm{\bm P_{\bm U_2} \bm \Delta {(\bm X^\dagger)}^T \bm \Delta \bm P_{\bm V_2}}_F^2 . \numberthis \label{equ:Delta3}
\end{align*}
Using the cyclic property of the trace and the idempotence property of $\bm P_{\bm V_2}$, the first term on the RHS of (\ref{equ:Delta3}) can be computed as
\begin{align*}
    &\norm{\bm X^\dagger \bm \Delta^T \bm P_{\bm U_2} \bm \Delta \bm P_{\bm V_2}}_F^2 = \tr\bigl(\bm X^\dagger \bm \Delta^T \bm P_{\bm U_2} \bm \Delta \bm P_{\bm V_2} \bm P_{\bm V_2} \bm \Delta^T \bm P_{\bm U_2} \bm \Delta (\bm X^\dagger)^T \bigr) = \tr \bigl( \bm \Delta^T \bm P_{\bm U_2} \bm \Delta \bm P_{\bm V_2} \bm \Delta^T \bm P_{\bm U_2} \bm \Delta (\bm X^\dagger)^T \bm X^\dagger \bigr) .
\end{align*}
Similarly, one can compute the second and the third terms on the RHS of (\ref{equ:Delta3}), then taking the expectation to obtain
\begin{align*}
    \E \bigl[ &\norm{\bm X^\dagger \bm \Delta^T \bm P_{\bm U_2} \bm \Delta \bm P_{\bm V_2} + \bm P_{\bm U_2} \bm \Delta \bm P_{\bm V_2} \bm \Delta^T \bm X^\dagger + \bm P_{\bm U_2} \bm \Delta {(\bm X^\dagger)}^T \bm \Delta \bm P_{\bm V_2}}_F^2 \bigr] \\
    &= \E \Bigl[ \tr \bigl( \bm \Delta^T \bm P_{\bm U_2} \bm \Delta \bm P_{\bm V_2} \bm \Delta^T \bm P_{\bm U_2} \bm \Delta (\bm X^\dagger)^T \bm X^\dagger \bigr) \Bigr] + \E \Bigl[ \tr \bigl( \bm \Delta^T \bm X^\dagger (\bm X^\dagger)^T \bm \Delta \bm P_{\bm V_2} \bm \Delta^T \bm P_{\bm U_2} \bm \Delta \bm P_{\bm V_2} \bigr) \Bigr] \\
    &\qquad \qquad \qquad + \E \Bigl[ \tr\bigl( \bm \Delta {(\bm X^\dagger)}^T \bm \Delta \bm P_{\bm V_2} \bm \Delta^T \bm X^\dagger \bm \Delta^T \bm P_{\bm U_2} \big) \Bigr]. \numberthis \label{equ:Delta3b}
\end{align*}
Next, to compute the three terms on the RHS of (\ref{equ:Delta3b}), we consider the following lemma:
\begin{lemma}\label{lem:ED4}
Assume the matrices $\bm A, \bm B, \bm C,$ and $\bm D$ in each of the following statements are of compatible dimensions such that the matrix product is valid. Then,
\begin{enumerate}
    \item $\E \bigl[\tr(\bm \Delta^T \bm A \bm \Delta \bm B \bm \Delta^T \bm C \bm \Delta \bm D) \bigr] = \tr(\bm A^T \bm C) \tr(\bm B^T \bm D) + \tr(\bm A \bm C) \tr(\bm B) \tr(\bm D) + \tr(\bm B \bm D) \tr(\bm A) \tr(\bm C))$,
    \item $\E \bigl[\tr(\bm \Delta \bm A \bm \Delta \bm B \bm \Delta^T \bm C \bm \Delta^T \bm D) \bigr] = \tr(\bm A^T \bm B \bm C^T \bm D) + \tr(\bm D \bm C \bm B \bm A) + \tr(\bm A \bm C) \tr(\bm B) \tr(\bm D)$.
\end{enumerate}
\end{lemma}
The proof of Lemma~\ref{lem:ED4} follows a similar derivation of the fourth-moment properties in \cite{neudecker1987fourth} and hence is omitted. Applying Lemma~\ref{lem:ED4} to the RHS of (\ref{equ:Delta3b}) and using the orthogonality between $\bm X^\dagger$ and $\bm P_{\bm U_2}, \bm P_{\bm V_2}$, we obtain
\begin{align*}
    \E \bigl[ &\norm{\bm X^\dagger \bm \Delta^T \bm P_{\bm U_2} \bm \Delta \bm P_{\bm V_2} + \bm P_{\bm U_2} \bm \Delta \bm P_{\bm V_2} \bm \Delta^T \bm X^\dagger + \bm P_{\bm U_2} \bm \Delta {(\bm X^\dagger)}^T \bm \Delta \bm P_{\bm V_2}}_F^2 \bigr] \\
    &= \tr(\bm P_{\bm U_2}) \tr(\bm P_{\bm V_2}) \tr((\bm X^\dagger)^T \bm X^\dagger) + \tr(\bm P_{\bm V_2}) \tr(\bm X^\dagger (\bm X^\dagger)^T) \tr(\bm P_{\bm U_2}) + \tr((\bm X^\dagger)^T \bm X^\dagger) \tr(\bm P_{\bm V_2}) \tr(\bm P_{\bm U_2}) \\
    &= 3\sigma^4 (m-r)(n-r)\norm{\bm X^\dagger}_F^2 , \numberthis \label{equ:M_hat2}
\end{align*}
where the last equality stems from $\tr(\bm P_{\bm U_2})=m-r$ and $\tr(\bm P_{\bm V_2})=n-r$.
Substituting (\ref{equ:M1}) and (\ref{equ:M_hat2}) into (\ref{equ:M2nd}) yields
\begin{align*}
    \E \Bigl[ \norm{\hat{\bm X}_2 - \bm X}_F^2 \Bigr] = \sigma^2 r(m+n-r) + 3\sigma^4 (m-r)(n-r)\norm{\bm X^\dagger}_F^2 . \numberthis \label{equ:M2}
\end{align*}

\item \textbf{The upper bound:}

From Corollary~\ref{COR:1ST}, we have
\begin{align*}
    \hat{\bm X} - \bm X = \P_r (\tilde{\bm X}) - \bm X = \bm \Delta - \bm P_{\bm U_2} \bm \Delta \bm P_{\bm V_2} + \bm R_{\bm X} (\bm \Delta) .
\end{align*}
Hence, by the triangle inequality, it holds that
\begin{align*}
    \norm{\hat{\bm X} - \bm X}_F \leq \norm{\bm \Delta - \bm P_{\bm U_2} \bm \Delta \bm P_{\bm V_2}}_F + \norm{\bm R_{\bm X} (\bm \Delta)}_F .
\end{align*}
Taking the expectation of the squared norm yields
\begin{align*}
    \E \Bigl[ \norm{\hat{\bm X} - \bm X}_F^2 \Bigr] &\leq \E \Bigl[ \Bigl( \norm{\bm \Delta - \bm P_{\bm U_2} \bm \Delta \bm P_{\bm V_2}}_F + \norm{\bm R_{\bm X} (\bm \Delta)}_F \Bigr)^2 \Bigr] . \numberthis \label{equ:M_bound}
\end{align*}


Applying Minkowski inequality \cite{hardy1952inequalities}, we can bound the RHS of (\ref{equ:M_hat2}) as
\begin{align*}
    \E \Bigl[ \Bigl( \norm{\bm \Delta - \bm P_{\bm U_2} \bm \Delta \bm P_{\bm V_2}}_F &+ \norm{\bm R_{\bm X} (\bm \Delta)}_F \Bigr)^2 \Bigr] \leq \biggl( \sqrt{\E \Bigl[ \norm{\bm \Delta - \bm P_{\bm U_2} \bm \Delta \bm P_{\bm V_2}}_F^2 \Bigr]} + \sqrt{\E \Bigl[ \norm{\bm R_{\bm X} (\bm \Delta)}_F^2 \Bigr]} \biggr)^2 . \numberthis \label{equ:min_bound}
\end{align*}
From (\ref{equ:boundR_2}), we can bound $\E \Bigl[ \norm{\bm R_{\bm X} (\bm \Delta)}_F^2 \Bigr]$ by
\begin{align*}
    \E \Bigl[ \norm{\bm R_{\bm X} (\bm \Delta)}_F^2 \Bigr] &\leq \E \Bigl[ \bigl(2(1+\sqrt{2}) \min \bigl\{ \frac{2}{\sigma_r} \norm{\bm \Delta}_F^2 , \norm{\bm \Delta}_F \bigr\} \bigr)^2 \Bigr] \\
    &\leq \min \biggl\{  \Bigl(4(1+\sqrt{2})\frac{\sigma}{\sigma_r}\Bigr)^2 \E \Bigl[ \norm{\bm \Delta}_F^4 \Bigr] , \Bigl(2(1+\sqrt{2})\Bigr)^2 \E \Bigl[ \norm{\bm \Delta}_F^2 \Bigr] \biggr\} , \numberthis \label{equ:R_min}
\end{align*}
where the last inequality is a special case of Jensen's inequality \cite{hardy1952inequalities} with the minimum of two linear functions as a concave function.
The fourth-order term on the RHS of (\ref{equ:R_min}) can be computed as
\begin{align*}
    \E \Bigl[ \norm{\bm \Delta}_F^4 \Bigr] = \E \Bigl[ \bigl( \sum_{i=1}^m \sum_{j=1}^n \Delta_{ij}^2 \bigr)^2 \Bigr] = \sum_{i,j,k,l} \E \Bigl[ \Delta_{ij}^2 \Delta_{kl}^2 \Bigr] = \sigma^4 \sum_{i,j,k,l} (1+2\delta_{ik}\delta_{jl}) = \sigma^4 (m^2 n^2 + 2mn) . \numberthis \label{equ:D4}
\end{align*}
Substituting (\ref{equ:M1}) and (\ref{equ:D4}) back into (\ref{equ:R_min}), then taking the square root, we have
\begin{align*}
    \sqrt{\E \Bigl[ \norm{\bm R_{\bm X} (\bm \Delta)}_F^2 \Bigr]} \leq \min \biggl\{  4(1+\sqrt{2})\frac{\sigma}{\sigma_r} \sqrt{m^2 n^2 + 2mn} , 2(1+\sqrt{2}) \sqrt{mn} \biggr\} .
\end{align*}
Substituting the bound in the last inequality and the equality in (\ref{equ:M1}) into (\ref{equ:min_bound}), we obtain the upper bound as
\begin{align*}
    \E \Bigl[ \norm{\hat{\bm X} &- \bm X}_F^2 \Bigr] \leq \sigma^2 \Biggl( \sqrt{r (m+n-r)} + \min \biggl\{  4(1+\sqrt{2})\frac{\sigma}{\sigma_r} \sqrt{m^2 n^2 + 2mn} , 2(1+\sqrt{2}) \sqrt{mn} \biggr\} \Biggr)^2 . \numberthis \label{equ:M_hat}
\end{align*}
Due to the nature of the bound given in (\ref{equ:boundR_2}), the bound in (\ref{equ:M_hat}) is taken as the minimum between a component that is linear in the norm of $\bm \Delta$ and a component that is quadratic in the norm of $\bm \Delta$.

\end{enumerate}


\begin{remark}
Asymptotically as $\sigma \to 0$, {all the ratios of the first-order approximation (\ref{equ:M1}), the second-order approximation (\ref{equ:M2}), and the upper bound (\ref{equ:M_hat}) to the MSE of the noisy matrix (\ref{equ:M_tilde}) converge to $r(m+n-r)/mn$. In general, this ratio is less than or equal to $1$, however, in low-rank scenarios it can be significantly smaller.} This indicates the TSVD estimator is effective in noise reduction when the noise is small, especially when the matrix $\bm X$ has low rank.
\end{remark}

\begin{remark}
The upper bound in (\ref{equ:M_hat}) attains the same value of the baseline $\sigma^2 mn$ when $\sigma=\sigma_2$, where
\begin{align*}
    \sigma_2 = \frac{\sigma_r\bigl(\sqrt{mn}-\sqrt{r(m+n-r)}\bigr)}{4(1+\sqrt{2})\sqrt{m^2 n^2 + 2mn}} , \numberthis\label{equ:point2}
\end{align*}
guaranteeing the superiority of the upper bound over the baseline in the case $\sigma < \sigma_2$.
\end{remark}

\begin{remark}
Let us define the $\rho$-knee point between two increasing functions of $\sigma$, e.g., $f(\sigma)$ and $g(\sigma)$, as the point at which $f(\sigma) = \rho g(\sigma)$ (for $\rho>1$).
Then, the $\rho$-knee point between the upper bound (\ref{equ:M_hat}) and the first-order approximation (\ref{equ:M1}) can be determined by
\begin{align*}
    \sigma_1 = \frac{\sigma_r(\sqrt{\rho}-1)\sqrt{r(m+n-r)}}{4(1+\sqrt{2})\sqrt{m^2 n^2 + 2mn}} . \numberthis\label{equ:point1}
\end{align*}
In addition, the $\rho$-knee point between the second-order approximation (\ref{equ:M2}) and the first-order approximation (\ref{equ:M1}) is given by
\begin{align*}
    \sigma_3 = \sqrt{\frac{(\rho-1)r(m+n-r)}{3(m-r)(n-r)\norm{\bm X^\dagger}_F^2}} > \sigma_1 . \numberthis\label{equ:point3}
\end{align*}
\end{remark}

\begin{figure}[t]
    \centering
    \includegraphics[width=0.6\textwidth]{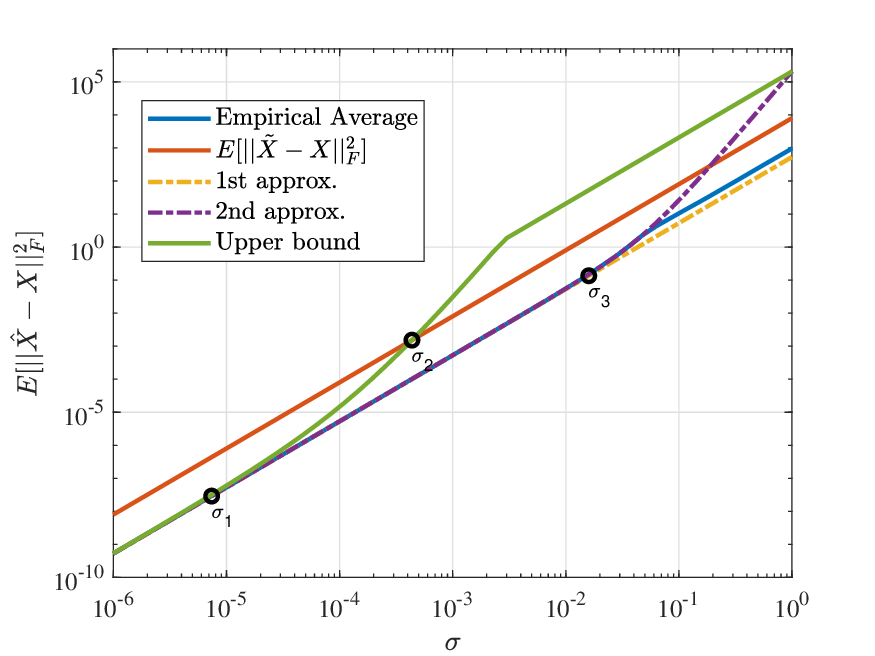}
    \caption{The MSE of the TSVD-based estimator $\hat{\bm X}$ for matrix denoising as a function of $\sigma$. The solid blue line represents the empirical estimate of MSE of $\hat{\bm X}$, i.e., $\E[ \norm{\hat{\bm X} - \bm X}_F^2 ]$. The solid red line is the MSE of the noisy matrix, i.e., $\E[\norm{\tilde{\bm X} - \bm X}_F^2] = \sigma^2 mn$. The dash-dotted yellow line and the dash-dotted purple line represent the first-order and second-order approximations of $\E[ \norm{\hat{\bm X} - \bm X}_F^2 ]$, i.e., $\E[ \norm{\hat{\bm X}_1 - \bm X}_F^2 ]$ and $\E[ \norm{\hat{\bm X}_2 - \bm X}_F^2 ]$, respectively. The solid green line is the upper bound on $\E[ \norm{\hat{\bm X} - \bm X}_F^2 ]$ given in (\ref{equ:M_hat}). The knee-points $\sigma_1$ and $\sigma_3$ represent the value of $\sigma$ for which the upper-bound and the second-order approximation deviate from the first order approximation by more than $10\%$, obtained by (\ref{equ:point1}) and (\ref{equ:point3}) with $\rho =1.1$. The point $\sigma_2$ is the intersection between the upper bound and the MSE of the noisy matrix, given by (\ref{equ:point2}).}
    \label{fig:denoising}
\end{figure}

Figure~\ref{fig:denoising} demonstrates the aforementioned analysis on the performance of the TSVD-based estimator for matrix denoising through a numerical experiment. 
{\bf Data generation.} We generate a matrix $\bm X$ with $m=100, n=80$, and $r=3$ by (i) taking the product of two random matrices, whose entries are $i.i.d.$ normally distributed ${\cal N}(0,1)$, of sizes $100 \times 3$ and $3 \times 80$, respectively; (ii) and dividing each entry of the obtained matrix by its Frobenius norm such that the resulting matrix satisfies $\norm{\bm X}_F=1$. In the experiment, we consider $51$ values of $\sigma$ in the interval of $[10^{-6},10^0]$, namely $\sigma \in \{10^{-6}, 10^{-5.88},10^{-5.76}, \ldots, 10^{0}\}$. 
For each value of $\sigma$, we compute the following quantities:
\begin{enumerate}
    \item the empirical MSE of the TSVD-based estimator $\E[\norm{\hat{\bm X} - \bm X}_F^2]$ by averaging the quantity $\norm{{\cal P}_r(\bm X + \bm \Delta) - \bm X}_F^2$ over $1000$ $i.i.d.$ instances of $\bm \Delta$,
    \item the MSE of the noisy matrix given in (\ref{equ:M_tilde}),
    \item the first-order approximation of the MSE of the TSVD-based estimator given in (\ref{equ:M1}),
    \item the second-order approximation of the MSE of the TSVD-based estimator given in (\ref{equ:M2nd}),
    \item the upper bound on the MSE of the TSVD-based estimator given in (\ref{equ:M_hat}).
\end{enumerate}
\noindent We display each of the aforementioned quantities as a function of $\sigma$ in Fig.~\ref{fig:denoising}. In addition, we calculate the points corresponding to $\sigma_1, \sigma_2$ and $\sigma_3$ using (\ref{equ:point1}), (\ref{equ:point2}), and (\ref{equ:point3}), respectively, with $\rho=1.1$, and include them in Fig.~\ref{fig:denoising}. 
\noindent {\bf Results and Analysis.} It can be observed from the plot that the empirical MSE of the TSVD-based estimator (solid blue) increases quadratically as a function of $\sigma$ (in the log-log scale, it appears as a straight line with slope equal to $2$).
The first-order approximation (dash-dotted yellow) and the second-order approximation (dash-dotted purple) match the empirical average well for $\sigma < \sigma_3 \approx 10^{-2}$.
In this range of $\sigma$, all of the three aforementioned quantities are lower than the MSE of the noisy matrix (solid red).
On the other hand, the upper bound (solid green) holds tightly when $\sigma < \sigma_1 \approx 10^{-4}$, providing an efficient guarantee on the performance of the TSVD-based estimator for denoising with the presence of small additive noises. However, as the noise variance increases, the upper bound appears loose, exceeding the MSE of the noisy matrix when $\sigma > \sigma_2 \approx 4 \times 10^{-4}$. The bound is developed for the worst-case noise scenario, in which the noise is adversarially selected to yield the largest perturbation error (see the proof of Theorem~\ref{THEO:RANKOP} in Appendix~\ref{apdx:rankop}) and not for the random noise case. Consequently, it is far more conservative, predicting a larger MSE than the actual MSE of the TSVD-based estimator.
Developing bounds for average-case scenario is a potential direction for future research.

\section{Conclusion}

In this paper, we derived a first-order perturbation expansion for the singular value truncation. When the underlying matrix has exact rank-$r$, we showed that the first-order approximation can be greatly simplified and further introduced a simple expression of the second-order perturbation expansion for the $r$-TSVD.
Next, we proposed an error bound on the first-order approximation of the $r$-TSVD about a rank-$r$ matrix. Our bound is universal in the sense that it holds for perturbation matrices with an arbitrary norm. Two open questions raised by our analysis are: \hlnew{(i) when the underlying matrix has arbitrary rank, whether there exists an explicit expression for the second-order perturbation expansion of the TSVD; (ii) and given the result in} Theorem~\ref{THEO:1ST}\hlnew{, whether it is possible to establish a global error bound on the first-order approximation of the $r$-TSVD.}

\section*{Declaration of competing interest}

There is no competing interest.

\section*{Acknowledgments}

This work is partially supported by the National Science Foundation grant CCF-1254218.

\appendix

\section{Auxiliary lemmas}
\label{apdx:aux}

This section summarizes some trivial results that will be used regularly in our subsequent derivation. \hlnew{The proofs of Lemmas~}\ref{lem:trivial}-\ref{lem:semi} can be found in \cite{meyer2000matrix} - Chapter~5. The proof of Lemma~\ref{lem:kron} can be found in \cite{graham2018kronecker} - Chapter~2.

\begin{lemma} \label{lem:trivial}
Assume the same setting as in Definition~\ref{def:SVD}. The following statements hold:
\begin{enumerate}
    \item $\bm U_1^T \bm U_1 = \bm V_1^T \bm V_1= \bm I_r$, $\bm U_2^T \bm U_2 = \bm I_{m-r}$, and $\bm V_2^T \bm V_2 = \bm I_{n-r}$, 
    \item $\bm U_1^T \bm U_2 = \bm 0_{r \times (m-r)}$ and $\bm V_1^T \bm V_2 = \bm 0_{r \times (n-r)}$, 
    \item $\bm P_{\bm U_1} \bm P_{\bm U_2} = \bm 0$ and $\bm P_{\bm V_1} \bm P_{\bm V_2} = \bm 0$.
\end{enumerate}
Furthermore, if $\bm X$ has rank $r$, then
\begin{enumerate}
    \item $\bm P_{\bm U_2} \bm X = \bm 0$ and $\bm X \bm P_{\bm V_2} = \bm 0$,
    \item $\bm X = \bm P_{\bm U_1} \bm X = \bm X \bm P_{\bm V_1}$,
    \item $\bm X (\bm X^\dagger)^T = \bm P_{\bm U_1}$ and $\bm X^T \bm X^\dagger = \bm P_{\bm V_1}$.
\end{enumerate}
\end{lemma}

\begin{lemma} \label{lem:Pr}
Assume the same setting as in Definition~\ref{def:SVD}. The following statements hold:
\begin{enumerate}
    \item $\P_r(\bm X) = \bm U_{1} \bm \Sigma_{1} \bm V_{1}^T = \bm P_{\bm U_{1}} \bm X = \bm X \bm P_{\bm V_{1}} = \bm P_{\bm U_{1}} \bm X \bm P_{\bm V_{1}}$,
    \item $\bm X - \P_r(\bm X) = \bm U_{2} \bm \Sigma_{2} \bm V_{2}^T = \bm P_{\bm U_{2}} \bm X = \bm X \bm P_{\bm V_{2}} = \bm P_{\bm U_{2}} \bm X \bm P_{\bm V_{2}}$.
\end{enumerate}
\end{lemma}


\begin{lemma} \label{lem:norm}
For any matrices $\bm A$ and $\bm B$ with compatible dimensions, the following inequalities hold
\begin{align*}
    \norm{\bm A \bm B}_2 \leq \norm{\bm A \bm B}_F \leq \min\{ \norm{\bm A}_F \norm{\bm B}_2, \norm{\bm A}_2 \norm{\bm B}_F \} \leq \norm{\bm A}_F \norm{\bm B}_F .
\end{align*}
\end{lemma}

\begin{lemma} \label{lem:Pythagoras}
(Pythagoras theorem for Frobenius norm) For any matrices $\bm A$ and $\bm B$ such that $\tr (\bm A^T \bm B) = 0$, it holds that
\begin{align*}
    \norm{\bm A + \bm B}_F = \sqrt{\norm{\bm A}_F^2 + \norm{\bm B}_F^2} .
\end{align*}
The matrices $\bm A$ and $\bm B$ in this case are said to be orthogonal to each other.
\end{lemma}

\begin{lemma} \label{lem:semi}
Let $\bm U$ be a semi-orthogonal matrix with orthonormal columns and $\bm P_{\bm U} = \bm U \bm U^T$. Then, for any matrices $\bm A$ and $\bm B$ that have compatible dimensions with $\bm U$, the followings hold
\begin{enumerate}
    \item $\norm{\bm U \bm A}_2 = \norm{\bm A}_2$ and $\norm{\bm U \bm A}_F = \norm{\bm A}_F$,
    \item $\norm{\bm B \bm U}_2 = \norm{\bm B \bm P_{\bm U}}_2 \leq \norm{\bm B}_2$ and $\norm{\bm B \bm U}_F = \norm{\bm B \bm P_{\bm U}}_F \leq \norm{\bm B}_F$.
\end{enumerate}
\end{lemma}

\begin{lemma} \label{lem:kron}
For any matrices $\bm A$, $\bm B$, $\bm C$, and $\bm D$ with compatible dimensions such that the matrix products are valid, the following holds
\begin{enumerate}
    \item $(\bm A \otimes \bm B) (\bm C \otimes \bm D) = (\bm A \bm C) \otimes (\bm B \bm D)$,
    \item $\vect(\bm A \bm B \bm C) = (\bm C^T \otimes \bm A) \vect(\bm B)$,
    \item $\norm{\bm A \otimes \bm B}_F = \norm{\bm A}_F \norm{\bm B}_F$,
    \item $\tr(\bm A \otimes \bm B) = \tr(\bm A) \tr(\bm B)$.
\end{enumerate}
\end{lemma}

\section{Proof of Theorem~\ref{THEO:1ST}}
\label{apdx:1st}

\hlnew{Recall that in this proof, we consider a matrix $\bm X$ having rank greater than or equal to $r$. With a slight abuse of notation, let us define $\bm R_{\bm X}(\bm \Delta)$ as follows:}
\begin{align}
    \bm R_{\bm X}(\bm \Delta) &= \P_r(\bm X + \bm \Delta) - \bigl( \P_r(\bm X) + \bm \Delta - \bm P_{\bm U_{2}} \bm \Delta \bm P_{\bm V_{2}} \bigr) \label{equ:G0} \\
    &= \bigl( \P_r(\bm X + \bm \Delta) - (\bm X + \bm \Delta) \bigr) + \bigl( \bm X - \P_r(\bm X) \bigr) + \bm P_{\bm U_{2}} \bm \Delta \bm P_{\bm V_{2}} . \label{equ:G1}
\end{align}
Since $\tilde{\bm X} = \bm X + \bm \Delta$, applying Lemma~\ref{lem:Pr} to (\ref{equ:G1}) yields
\begin{align*}
    \bm R_{\bm X}(\bm \Delta) &= -\bm P_{\tilde{\bm U}_2} \tilde{\bm X} \bm P_{\tilde{\bm V}_2} + \bm P_{\bm U_2} \bm X \bm P_{\bm V_2} + \bm P_{\bm U_{2}} \bm \Delta \bm P_{\bm V_{2}} \\
    & = -\bm P_{\tilde{\bm U}_2} \tilde{\bm X} \bm P_{\tilde{\bm V}_2} + \bm P_{\bm U_2} \tilde{\bm X} \bm P_{\bm V_2} .
\end{align*}
Denote $\bm \delta_{\bm P_{\bm U_2}} = \bm P_{\tilde{\bm U}_2} - \bm P_{\bm U_2}$ and $\bm \delta_{\bm P_{\bm V_2}} = \bm P_{\tilde{\bm V}_2} - \bm P_{\bm U_2}$. By rewriting $\bm P_{\tilde{\bm U}_2} = \bm P_{\bm U_2} + \bm \delta_{\bm P_{\bm U_2}}$ and $\bm P_{\tilde{\bm V}_2} = \bm P_{\bm V_2} + \bm \delta_{\bm P_{\bm V_2}}$, we can further simplify the last equation as
\begin{align} \label{equ:G2}
    \bm R_{\bm X}(\bm \Delta) = -\bm \delta_{\bm P_{\bm U_2}} \tilde{\bm X} \bm P_{\bm V_2} - \bm P_{\bm U_2} \tilde{\bm X} \bm \delta_{\bm P_{\bm V_2}} - \bm \delta_{\bm P_{\bm U_2}} \tilde{\bm X} \bm \delta_{\bm P_{\bm V_2}} .
\end{align}

\begin{lemma} \label{lem:dP}
The perturbations of singular subspaces satisfy
\begin{subequations}
\label{equ:dPUV}%
\begin{align}
    &\bm \delta_{\bm P_{\bm U_2}} = \bm U_1 \bm Q^T \bm U_2^T + \bm U_2 \bm Q \bm U_1^T + \bm \O(\norm{\bm \Delta}_F^2) , \label{equ:dPU} \\
    &\bm \delta_{\bm P_{\bm V_2}} = - \bm V_1 \bm P^T \bm V_2^T - \bm V_2 \bm P \bm V_1^T + \bm \O(\norm{\bm \Delta}_F^2) . \label{equ:dPV} 
\end{align}
\end{subequations}
\end{lemma}
\noindent The proof of Lemma~\ref{lem:dP} is given at the end of this section. From this lemma, it is clear that $\bm \delta_{\bm P_{\bm U_2}}$ and $\bm \delta_{\bm P_{\bm V_2}}$ are in the order of $\norm{\bm \Delta}_F$. Substituting $\tilde{\bm X} = \bm X + \bm \Delta$ into (\ref{equ:G2}) and collecting second-order terms yield
\begin{align} \label{equ:G3}
    \bm R_{\bm X}(\bm \Delta) = -\bm \delta_{\bm P_{\bm U_2}} \bm X \bm P_{\bm V_2} - \bm P_{\bm U_2} \bm X \bm \delta_{\bm P_{\bm V_2}} + \bm \O(\norm{\bm \Delta}_F^2) .
\end{align}
Substituting (\ref{equ:dPU}) into the first term on the RHS of (\ref{equ:G3}), we obtain
\begin{align*}
    \bm \delta_{\bm P_{\bm U_2}} \bm X \bm P_{\bm V_2} &= \bigl( \bm U_1 \bm Q^T \bm U_2^T + \bm U_2 \bm Q \bm U_1^T \bigr) \bm U_2 \bm \Sigma_2 \bm V_2^T + \bm \O(\norm{\bm \Delta}_F^2) .
\end{align*}
Since $\bm U_2^T \bm U_2 = \bm I$ and $\bm U_1^T \bm U_2 = \bm 0$, we further have
\begin{align} \label{equ:G4}
    \bm \delta_{\bm P_{\bm U_2}} \bm X \bm P_{\bm V_2} &= \bm U_1 \bm Q^T \bm \Sigma_2 \bm V_2^T + \bm \O(\norm{\bm \Delta}_F^2) .
\end{align}
Similarly, the second term on the RHS of (\ref{equ:G3}) can be represented as
\begin{align} \label{equ:G5}
    \bm P_{\bm U_2} \bm X \bm \delta_{\bm P_{\bm V_2}} &= -\bm U_2 \bm \Sigma_2 \bm P \bm V_1^T + \bm \O(\norm{\bm \Delta}_F^2) .
\end{align}
Substituting (\ref{equ:G4}) and (\ref{equ:G5}) back into (\ref{equ:G3}), we have
\begin{align} \label{equ:G6}
    \bm R_{\bm X}(\bm \Delta) = - \bm U_1 \bm Q^T \bm \Sigma_2 \bm V_2^T + \bm U_2 \bm \Sigma_2 \bm P \bm V_1^T + \bm \O(\norm{\bm \Delta}_F^2) .
\end{align}
Now we can vectorize (\ref{equ:G6}) and apply Lemma~\ref{lem:kron} to obtain
\begin{align} \label{equ:G7}
    \vect \bigl( \bm R_{\bm X}(\bm \Delta) \bigr) = (\bm V_2 \bm \Sigma_2^T \otimes \bm U_1) \vect(-\bm Q^T) + (\bm V_1 \otimes \bm U_2 \bm \Sigma_2) \vect(\bm P) + \bm \O(\norm{\bm \Delta}_F^2) .
\end{align}
Let us now consider each term on the RHS of (\ref{equ:G7}). From Proposition~\ref{pro:QP}, it follows that
\begin{align} \label{equ:G8}
    \vect(-\bm Q^T) &= (\bm I_{m-r} \otimes \bm \Sigma_1^2 - \bm \Sigma_2 \bm \Sigma_2^T \otimes \bm I_r)^{-1} \vect(\bm E_{12} \bm \Sigma_2^T + \bm \Sigma_1 \bm E_{21}^T) + \bm \O(\norm{\bm \Delta}_F^2) .
\end{align}
Replacing $\bm E_{ij} = \bm U_i^T \bm \Delta \bm V_j$, for $i,j \in \{1,2\}$, and using Lemma~\ref{lem:kron}, (\ref{equ:G8}) becomes
\begin{align} \label{equ:G9}
    \vect(-\bm Q^T) = (\bm I_{m-r} \otimes \bm \Sigma_1^2 - \bm \Sigma_2 \bm \Sigma_2^T \otimes \bm I_r)^{-1} \bigl( (\bm \Sigma_2 \bm V_2^T \otimes \bm U_1^T) \vect(\bm \Delta) + (\bm U_2^T \otimes \bm \Sigma_1 \bm V_1^T) \vect(\bm \Delta^T) \bigr) + \bm \O(\norm{\bm \Delta}_F^2) .
\end{align}
Since $\bm \Sigma_1$ and $\bm \Sigma_2$ are diagonal, so is $(\bm I_{m-r} \otimes \bm \Sigma_1^2 - \bm \Sigma_2 \bm \Sigma_2^T \otimes \bm I_r)^{-1}$. The following lemma provides an insight into the structure of this inversion.

\begin{lemma} \label{lem:D}
Let $\bm D = (\bm I_{m-r} \otimes \bm \Sigma_1^2 - \bm \Sigma_2 \bm \Sigma_2^T \otimes \bm I_r)^{-1}$. Then
\begin{align*}
    \bm D = \sum_{i=1}^r \sum_{k=1}^{m-r} d_{ik} \bigl(\bm e^{m-r}_{k} (\bm e^{m-r}_{k})^T \bigr) \otimes \bigl(\bm e^r_{i} (\bm e^r_{i})^T \bigr) ,
\end{align*}
where $d_{ik} = \frac{1}{\sigma_i^2 - \sigma_{r+k}^2}$, for $i=1,\ldots,r$ and $k=1,\ldots,m-r$. 
\end{lemma}
\noindent The proof of Lemma~\ref{lem:D} is given at the end of this section. Now using Lemma~\ref{lem:D} and left-multiplying both sides of (\ref{equ:G9}) by $(\bm V_2 \bm \Sigma_2^T \otimes \bm U_1)$, we obtain
\begin{align*}
    (\bm V_2 \bm \Sigma_2^T &\otimes \bm U_1) \vect(-\bm Q^T) = \sum_{i=1}^r \sum_{k=1}^{m-r} d_{ik} (\bm V_2 \bm \Sigma_2^T \otimes \bm U_1) \Bigl( \bigl(\bm e^{m-r}_{k} (\bm e^{m-r}_{k})^T \bigr) \otimes \bigl(\bm e^r_{i} (\bm e^r_{i})^T \bigr) \Bigr) (\bm \Sigma_2 \bm V_2^T \otimes \bm U_1^T) \vect(\bm \Delta) \\
    &+ \sum_{i=1}^r \sum_{k=1}^{m-r} d_{ik} (\bm V_2 \bm \Sigma_2^T \otimes \bm U_1) \Bigl( \bigl(\bm e^{m-r}_{k} (\bm e^{m-r}_{k})^T \bigr) \otimes \bigl(\bm e^r_{i} (\bm e^r_{i})^T \bigr) \Bigr) (\bm U_2^T \otimes \bm \Sigma_1 \bm V_1^T) \vect(\bm \Delta^T) + \bm \O(\norm{\bm \Delta}_F^2) . \numberthis \label{equ:G10}
\end{align*}
Moreover, applying Lemma~\ref{lem:kron}-1, we have 
\begin{align*}
    (\bm V_2 \bm \Sigma_2^T \otimes \bm U_1) \Bigl( \bigl(\bm e^{m-r}_{k} (\bm e^{m-r}_{k})^T \bigr) \otimes \bigl(\bm e^r_{i} (\bm e^r_{i})^T \bigr) \Bigr) (\bm \Sigma_2 \bm V_2^T \otimes \bm U_1^T) &= \bigl( \bm V_2 \bm \Sigma_2^T \bm e^{m-r}_{k} (\bm e^{m-r}_{k})^T \bm \Sigma_2 \bm V_2^T \bigr) \otimes \bigl( \bm U_1 \bm e^r_{i} (\bm e^r_{i})^T \bm U_1^T \bigr) \\
    &= \sigma_{r+k}^2 (\bm v_{r+k} \bm v_{r+k}^T) \otimes (\bm u_i \bm u_i^T) , \numberthis \label{equ:G11} 
\end{align*}
and similarly,
\begin{align*}
    (\bm V_2 \bm \Sigma_2^T \otimes \bm U_1) \Bigl( \bigl(\bm e^{m-r}_{k} (\bm e^{m-r}_{k})^T \bigr) \otimes \bigl(\bm e^r_{i} (\bm e^r_{i})^T \bigr) \Bigr) (\bm U_2^T \otimes \bm \Sigma_1 \bm V_1^T) = \sigma_i \sigma_{r+k} (\bm v_{r+k} \bm u_{r+k}^T) \otimes (\bm u_i \bm v_i^T) . \numberthis \label{equ:G12}
\end{align*}
Substituting (\ref{equ:G11}) and (\ref{equ:G12}) back into (\ref{equ:G10}) and performing a change of variable $j=r+k$, we obtain
\begin{align*}
    (\bm V_2 \bm \Sigma_2^T \otimes \bm U_1) \vect(-\bm Q^T) = \sum_{i=1}^r &\sum_{j=r+1}^m \frac{\sigma_j^2}{\sigma_i^2 - \sigma_j^2} (\bm v_j \bm v_j^T) \otimes (\bm u_i \bm u_i^T) \vect(\bm \Delta) \\
    &+ \sum_{i=1}^r \sum_{j=r+1}^m \frac{\sigma_i \sigma_j}{\sigma_i^2 - \sigma_j^2} (\bm v_j \bm u_j^T) \otimes (\bm u_i \bm v_i^T) \vect(\bm \Delta^T) + \bm \O(\norm{\bm \Delta}_F^2) . \numberthis \label{equ:G13}
\end{align*}
Following a similar derivation, we also have
\begin{align*}
    (\bm V_1 \otimes \bm U_2 \bm \Sigma_2) \vect(\bm P) = \sum_{i=1}^r &\sum_{j=r+1}^m \frac{\sigma_j^2}{\sigma_i^2 - \sigma_j^2} (\bm v_i \bm v_i^T) \otimes (\bm u_j \bm u_j^T) \vect(\bm \Delta) \\
    &+ \sum_{i=1}^r \sum_{j=r+1}^m \frac{\sigma_i \sigma_j}{\sigma_i^2 - \sigma_j^2} (\bm v_i \bm u_i^T) \otimes (\bm u_j \bm v_j^T) \vect(\bm \Delta^T) + \bm \O(\norm{\bm \Delta}_F^2) . \numberthis \label{equ:G14}
\end{align*}
Substituting (\ref{equ:G13}) and (\ref{equ:G14}) back into (\ref{equ:G7}) yields
\begin{align*}
    \vect \bigl( \bm R_{\bm X}(\bm \Delta) \bigr) = \sum_{i=1}^r \sum_{j=r+1}^m &\biggl( \frac{\sigma_j^2}{\sigma_i^2 - \sigma_j^2} \bigl( (\bm v_j \bm v_j^T) \otimes (\bm u_i \bm u_i^T) + (\bm v_i \bm v_i^T) \otimes (\bm u_j \bm u_j^T) \bigr) \vect(\bm \Delta) \\
    &+ \frac{\sigma_i \sigma_j}{\sigma_i^2 - \sigma_j^2} \bigl( (\bm v_j \bm u_j^T) \otimes (\bm u_i \bm v_i^T) + (\bm v_i \bm u_i^T) \otimes (\bm u_j \bm v_j^T) \bigr) \vect(\bm \Delta^T) \biggr) + \bm \O(\norm{\bm \Delta}_F^2) . \numberthis \label{equ:G15}
\end{align*}
Truncating the inner summation, with $\sigma_j = 0$ for $j>n$, and applying Lemma~\ref{lem:kron}-2 to the RHS of (\ref{equ:G15}), we obtain
\begin{align*}
    \bm R_{\bm X}(\bm \Delta) = \sum_{i=1}^r \sum_{j=r+1}^n &\biggl( \frac{\sigma_{j}^2}{\sigma_{i}^2-\sigma_{j}^2} (\bm u_{i} \bm u_{i}^T \bm \Delta \bm v_{j} \bm v_{j}^T + \bm u_{j} \bm u_{j}^T \bm \Delta \bm v_{i} \bm v_{i}^T) + \frac{\sigma_{i} \sigma_{j}}{\sigma_{i}^2-\sigma_{j}^2} (\bm u_{i} \bm v_{i}^T \bm \Delta^T \bm u_{j} \bm v_{j}^T + \bm u_{j} \bm v_{j}^T \bm \Delta^T \bm u_{i} \bm v_{i}^T) \biggr) + \bm \O(\norm{\bm \Delta}_F^2) .
\end{align*}
Our theorem now follows on the definition of $\bm R_{\bm X}(\bm \Delta)$ in (\ref{equ:G0}).

\subsection{Proof of Lemma~\ref{lem:dP}}
Using the fact from Proposition~\ref{pro:stewart} that $\bm P_{\tilde{\bm U}_2} = \bm P_{\hat{\bm U}_2}$, we can re-express the subspace difference as 
\begin{align} \label{equ:dP1a}
    \bm \delta_{\bm P_{U_2}} = \bm P_{\tilde{\bm U}_2} - \bm P_{\bm U_2} = \bm P_{\hat{\bm U}_2} - \bm P_{\bm U_2} = \hat{\bm U}_2 \hat{\bm U}_2^T - \bm U_2 \bm U_2^T .
\end{align}
Substituting (\ref{equ:perturbed_U2}) into (\ref{equ:dP1a}) yields
\begin{align} \label{equ:dP1}
    \bm \delta_{\bm P_{U_2}} = (\bm U_2 + \bm U_1 \bm Q^T) (\bm I_{m-r} + \bm Q \bm Q^T)^{-1} (\bm U_2^T + \bm Q \bm U_1^T) - \bm U_2 \bm U_2^T .
\end{align}
Since $\bm Q = \bm \O(\norm{\bm \Delta}_F)$ and $(\bm I_{m-r} + \bm Q \bm Q^T)^{-1} = \bm I_{m-r} - \bm Q \bm Q^T (\bm I_{m-r} + \bm Q \bm Q^T)^{-1} = \bm I_{m-r} + \bm \O(\norm{\bm \Delta}_F^2)$, (\ref{equ:dP1}) can be simplified by absorbing second-order terms:
\begin{align*}
    \bm \delta_{\bm P_{U_2}} &= (\bm U_2 + \bm U_1 \bm Q^T) (\bm U_2^T + \bm Q \bm U_1^T) - \bm U_2 \bm U_2^T + \bm \O(\norm{\bm \Delta}_F^2) \\
    &= \bm U_1 \bm Q^T \bm U_2^T + \bm U_2 \bm Q \bm U_1^T + \bm U_1 \bm Q^T \bm Q \bm U_1^T + \bm \O(\norm{\bm \Delta}_F^2) \\
    &= \bm U_1 \bm Q^T \bm U_2^T + \bm U_2 \bm Q \bm U_1^T + \bm \O(\norm{\bm \Delta}_F^2) .
\end{align*}
The equation $\bm \delta_{\bm P_{\bm V_2}} = - \bm V_1 \bm P^T \bm V_2^T - \bm V_2 \bm P \bm V_1^T + \bm \O(\norm{\bm \Delta}_F^2)$ can be proved by a similar derivation. Since $\bm Q$ and $\bm P$ are in the order of $\norm{\bm \Delta}_F$, so do $\bm \delta_{\bm P_{\bm U_2}}$ and $\bm \delta_{\bm P_{\bm V_2}}$.

\subsection{Proof of Lemma~\ref{lem:D}}
Recall that 
\begin{align*}
    \bm \Sigma_1^2 = \begin{bmatrix*}[r] \sigma_1^2 & \ldots & 0 \\ & \ddots & \\ 0 & \ldots & \sigma_r^2 \end{bmatrix*} \in \R^{r \times r} \qquad \text{and} \qquad \bm \Sigma_2 \bm \Sigma_2^T = \begin{bmatrix*}[r] \sigma_{r+1}^2 & \ldots & 0 \\ & \ddots & \\ 0 & \ldots & \sigma_m^2 \end{bmatrix*} \in \R^{(m-r) \times (m-r)} .
\end{align*}
By the definition of the Kronecker product, we have
\begin{align*}
    \bm I_{m-r} \otimes \bm \Sigma_1^2 - \bm \Sigma_2 \bm \Sigma_2^T \otimes \bm I_r = \begin{bmatrix*} \bm \Sigma_1^2 - \sigma_{r+1}^2 \bm I_r & \ldots & \bm 0_r \\ & \ddots & \\ \bm 0_r & \ldots & \bm \Sigma_1^2 - \sigma_{m}^2 \bm I_r \end{bmatrix*} \in \R^{(m-r)r \times (m-r)r} .
\end{align*}
Therefore, we can invert this diagonal matrix by considering each of the $r \times r$ blocks:
\begin{align*}
    \bm D &= (\bm I_{m-r} \otimes \bm \Sigma_1^2 - \bm \Sigma_2 \bm \Sigma_2^T \otimes \bm I_r)^{-1} \\
    &= \begin{bmatrix*} (\bm \Sigma_1^2 - \sigma_{r+1}^2 \bm I_r)^{-1} & \ldots & \bm 0_r \\ & \ddots & \\ \bm 0_r & \ldots & (\bm \Sigma_1^2 - \sigma_{m}^2 \bm I_r)^{-1} \end{bmatrix*} .
\end{align*}
Now it is easy to verify that, for $i=1,\ldots,r$ and $k=1,\ldots,m-r$, the $i$-th diagonal entry of the $k$-th diagonal block, is $d_{ik} = 1/(\sigma_i^2 - \sigma_{r+k}^2)$. Furthermore, since $\bigl(\bm e^{m-r}_{k} (\bm e^{m-r}_{k})^T \bigr) \otimes \bigl(\bm e^r_{i} (\bm e^r_{i})^T \bigr)$ is a $(m-r)r \times (m-r)r$ matrix of all zeros but the $i$-th diagonal entry of the $k$-th diagonal block is $1$, we represent $\bm D$ as the sum of $(m-r)r$ rank-$1$ matrices:
\begin{align*}
    \bm D = \sum_{i=1}^r \sum_{k=1}^{m-r} d_{ik} \bigl(\bm e^{m-r}_{k} (\bm e^{m-r}_{k})^T \bigr) \otimes \bigl(\bm e^r_{i} (\bm e^r_{i})^T \bigr) .
\end{align*}

\section{Proof of Theorem~\ref{THEO:2ND}}
\label{apdx:2nd}

By the definition of the $r$-TSVD in (\ref{equ:Eckart}), we have
\begin{align} \label{equ:Pr10}
    \P_r (\tilde{\bm X}) = \bm P_{\tilde{\bm U}_1} \tilde{\bm X} \bm P_{\tilde{\bm V}_1} .
\end{align}
Since we assume $\bm X$ has exact rank $r$, the perturbed matrix can be represented as $\tilde{\bm X} = \bm X + \bm \Delta = \bm U_1 \bm \Sigma_1 \bm V_1^T + \bm \Delta$.
\hlnew{Substituting this back} into (\ref{equ:Pr10}) yields
\begin{align} \label{equ:Pr1} 
    \P_r(\bm X + \bm \Delta) = \bm P_{\tilde{\bm U}_1} (\bm U_1 \bm \Sigma_1 \bm V_1^T + \bm \Delta) \bm P_{\tilde{\bm V}_1} .
\end{align}
Similar to the derivation of (\ref{equ:dP1}), we obtain
$\bm P_{\tilde{\bm U}_1} = (\bm U_1 - \bm U_2 \bm Q)(\bm I_r + \bm Q^T \bm Q)^{-1} (\bm U_1^T - \bm Q^T \bm U_2^T)$ and $\bm P_{\tilde{\bm V}_1} = (\bm V_1 + \bm V_2\bm P) (\bm I_r + \bm P^T \bm P)^{-1} (\bm V_1^T + \bm P^T \bm V_2^T)$.
Substituting the expressions of $\bm P_{\tilde{\bm U}_1}$ and $\bm P_{\tilde{\bm V}_1}$ back into (\ref{equ:Pr1}), we obtain
\begin{align*}
    \P_r(\bm X + \bm \Delta) = (\bm U_1 - \bm U_2 \bm Q) (\bm I_r + \bm Q^T \bm Q)^{-1} &(\bm U_1^T - \bm Q^T \bm U_2^T) (\bm U_1 \bm \Sigma_1 \bm V_1^T + \bm \Delta) \\
    &\cdot (\bm V_1 + \bm V_2\bm P) (\bm I_r + \bm P^T \bm P)^{-1} (\bm V_1^T + \bm P^T \bm V_2^T) . \numberthis \label{equ:Pr11}
\end{align*}
By orthogonality, the product of three terms in the middle of the RHS of (\ref{equ:Pr11}) can be expanded and simplified as
\begin{align*}
    (\bm U_1^T - \bm Q^T \bm U_2^T) (\bm U_1 \bm \Sigma_1 \bm V_1^T + \bm \Delta) (\bm V_1 + \bm V_2\bm P) = (\bm \Sigma_1 + \bm E_{11}) + (\bm E_{12} \bm P - \bm Q^T \bm E_{21} - \bm Q^T \bm E_{22} \bm P) .
\end{align*}
Therefore, (\ref{equ:Pr11}) is equivalent to
\begin{align} \label{equ:Pr2}
    \P_r(\bm X + \bm \Delta) &= (\bm U_1 - \bm U_2 \bm Q) (\bm I_r + \bm Q^T \bm Q)^{-1} (\bm \Sigma_1 + \bm E_{11}) (\bm I_r + \bm P^T \bm P)^{-1} (\bm V_1^T + \bm P^T \bm V_2^T) \nonumber \\
    &\qquad + (\bm U_1 - \bm U_2 \bm Q) (\bm I_r + \bm Q^T \bm Q)^{-1} (\bm E_{12} \bm P - \bm Q^T \bm E_{21} - \bm Q^T \bm E_{22} \bm P) (\bm I_r + \bm P^T \bm P)^{-1} (\bm V_1^T + \bm P^T \bm V_2^T) .
\end{align}
Let us first focus on the first term on the RHS of (\ref{equ:Pr2}). 
Similar to the result after (\ref{equ:dP1}), we have $(\bm I_r + \bm Q^T \bm Q)^{-1} = \bm I_r - (\bm I_r + \bm Q^T \bm Q)^{-1} \bm Q^T \bm Q$ and $(\bm I_r + \bm P^T \bm P)^{-1} = \bm I_r - \bm P^T \bm P (\bm I_r + \bm P^T \bm P)^{-1}$, and hence
\begin{align*}
    (\bm U_1 &- \bm U_2 \bm Q) (\bm I_r + \bm Q^T \bm Q)^{-1} (\bm \Sigma_1 + \bm E_{11}) (\bm I_r + \bm P^T \bm P)^{-1} \\
    &= (\bm U_1 - \bm U_2 \bm Q) \Bigr(\bm I_r - (\bm I_r + \bm Q^T \bm Q)^{-1} \bm Q^T \bm Q \Bigr) (\bm \Sigma_1 + \bm E_{11}) \Bigl( \bm I_r - \bm P^T \bm P (\bm I_r + \bm P^T \bm P)^{-1} \Bigr) (\bm V_1^T + \bm P^T \bm V_2^T) \\
    &= (\bm U_1 - \bm U_2 \bm Q) (\bm \Sigma_1 + \bm E_{11}) (\bm V_1^T + \bm P^T \bm V_2^T) - (\bm U_1 - \bm U_2 \bm Q) (\bm \Sigma_1 + \bm E_{11}) \bm P^T \bm P (\bm I_r + \bm P^T \bm P)^{-1} (\bm V_1^T + \bm P^T \bm V_2^T) \\
    &\qquad -(\bm U_1 - \bm U_2 \bm Q) (\bm I_r + \bm Q^T \bm Q)^{-1} \bm Q^T \bm Q (\bm \Sigma_1 + \bm E_{11}) (\bm I_r + \bm P^T \bm P)^{-1} (\bm V_1^T + \bm P^T \bm V_2^T) . \numberthis \label{equ:Pr3}
\end{align*}
Recall that $\bm X = \bm U_1 \bm \Sigma_1 \bm V_1^T$ and $\bm E_{11} = \bm U_1^T \bm \Delta \bm V_1$. The product $(\bm U_1 - \bm U_2 \bm Q) (\bm \Sigma_1 + \bm E_{11}) (\bm V_1^T + \bm P^T \bm V_2^T)$ can be expanded as
\begin{align*}
    (\bm U_1 - \bm U_2 \bm Q) &(\bm \Sigma_1 + \bm E_{11}) (\bm V_1^T + \bm P^T \bm V_2^T) \\
    &= \bm X + \bm U_1 \bm E_{11} \bm V_1^T + \bm U_1 (\bm \Sigma_1 + \bm E_{11}) \bm P^T \bm V_2^T - \bm U_2\bm Q (\bm \Sigma_1 + \bm E_{11}) \bm V_1^T - \bm U_2 \bm Q (\bm \Sigma_1 + \bm E_{11}) \bm P^T \bm V_2^T . \numberthis \label{equ:PrX1}
\end{align*}
In order to make up the first-order terms that involve $\bm \Delta$, we need to decompose the perturbation into $4$ components corresponding to different subspaces as follows.
Since $\bm P_{\bm U_{1}} + \bm P_{\bm U_{2}} = \bm I_m$ and $\bm P_{\bm V_{1}} + \bm P_{\bm V_{2}} = \bm I_n$, we have
\begin{align*}
    \bm \Delta &= \bm P_{\bm U_1} \bm \Delta \bm P_{\bm V_1} + \bm P_{\bm U_2} \bm \Delta \bm P_{\bm V_1} + \bm P_{\bm U_1} \bm \Delta \bm P_{\bm V_2} + \bm P_{\bm U_2} \bm \Delta \bm P_{\bm V_2} . \numberthis \label{equ:PrX3}
\end{align*}
Reorganizing terms in (\ref{equ:PrX3}) as
\begin{align*}
    \bm P_{\bm U_1} \bm \Delta \bm P_{\bm V_1} = \bm \Delta - \bm P_{\bm U_2} \bm \Delta \bm P_{\bm V_2} - \bm P_{\bm U_1} \bm \Delta \bm P_{\bm V_2} - \bm P_{\bm U_2} \bm \Delta \bm P_{\bm V_1} ,
\end{align*}
and using the definition of $\bm E$ in (\ref{equ:E}), we further have
\begin{align*}
    \bm U_1 \bm E_{11} \bm V_1^T &= \bm \Delta - \bm P_{\bm U_2} \bm \Delta \bm P_{\bm V_2} - \bm U_1 \bm E_{12} \bm V_2^T - \bm U_2 \bm E_{21} \bm V_1^T . \numberthis \label{equ:PrX2}
\end{align*}
Thus, substituting (\ref{equ:PrX2}) back into (\ref{equ:PrX1}) and rearranging terms yield
\begin{align*}
    (\bm U_1 - \bm U_2 \bm Q) (\bm \Sigma_1 + \bm E_{11}) (\bm V_1^T + \bm P^T \bm V_2^T) 
    &= \bm X + \bm \Delta - \bm P_{\bm U_2} \bm \Delta \bm P_{\bm V_2} + \bm U_1 \bigl( (\bm \Sigma_1 + \bm E_{11}) \bm P^T - \bm E_{12} \bigr) \bm V_2^T \\
    &\qquad - \bm U_2 \bigl(\bm Q (\bm \Sigma_1 + \bm E_{11}) + \bm E_{21} \bigr) \bm V_1^T - \bm U_2 \bm Q (\bm \Sigma_1 + \bm E_{11}) \bm P^T \bm V_2^T . \numberthis \label{equ:Pr4}
\end{align*}
Substituting (\ref{equ:Pr3}) and (\ref{equ:Pr4}) back into (\ref{equ:Pr2}), we obtain
\begin{align*}
    \P_r(\bm X + \bm \Delta) &= \bm X + \bm \Delta - \bm P_{\bm U_2} \bm \Delta \bm P_{\bm V_2} \\
    &\qquad + \bm U_1 \bigl( (\bm \Sigma_1 + \bm E_{11}) \bm P^T - \bm E_{12} \bigr) \bm V_2^T - \bm U_2 \bigl(\bm Q (\bm \Sigma_1 + \bm E_{11}) + \bm E_{21} \bigr) \bm V_1^T - \bm U_2 \bm Q (\bm \Sigma_1 + \bm E_{11}) \bm P^T \bm V_2^T \\
    &\qquad + (\bm U_1 - \bm U_2 \bm Q) (\bm I_r + \bm Q^T \bm Q)^{-1} \Bigl( - (\bm I_r + \bm Q^T \bm Q) (\bm \Sigma_1 + \bm E_{11}) \bm P^T \bm P + \bm Q^T \bm Q (\bm \Sigma_1 + \bm E_{11}) \\
    &\qquad \qquad + (\bm E_{12} \bm P - \bm Q^T \bm E_{21} - \bm Q^T \bm E_{22} \bm P) \Bigr) (\bm I_r + \bm P^T \bm P)^{-1} (\bm V_1^T + \bm P^T \bm V_2^T) . \numberthis \label{equ:Pr151}
\end{align*}
Applying (\ref{equ:QPS}), we have
\begin{align*}
    \bm U_1 \bigl( (\bm \Sigma_1 &+ \bm E_{11}) \bm P^T - \bm E_{12} \bigr) \bm V_2^T - \bm U_2 \bigl(\bm Q (\bm \Sigma_1 + \bm E_{11}) + \bm E_{21} \bigr) \bm V_1^T - \bm U_2 \bm Q (\bm \Sigma_1 + \bm E_{11}) \bm P^T \bm V_2^T \\
    &= \bm U_1 \bm Q^T (\bm E_{21} \bm P^T - \bm E_{22}) \bm V_2^T + \bm U_2 (\bm E_{22} - \bm Q \bm E_{12}) \bm P \bm V_1^T + \bm U_2 (\bm E_{21} + \bm E_{22} \bm P + \bm Q \bm E_{12} \bm P) \bm P^T \bm V_2^T , \numberthis \label{equ:Pr152}
\end{align*}
and
\begin{align*}    
    - (\bm I_r + \bm Q^T \bm Q) &(\bm \Sigma_1 + \bm E_{11}) \bm P^T \bm P + \bm Q^T \bm Q (\bm \Sigma_1 + \bm E_{11}) + (\bm E_{12} \bm P - \bm Q^T \bm E_{21} - \bm Q^T \bm E_{22} \bm P) \\
    &= \bigl( \bm E_{12} \bm P - (\bm \Sigma_1 + \bm E_{11}) \bm P^T \bm P \bigr) - \bigl( \bm Q^T \bm E_{21} + \bm Q^T \bm Q (\bm \Sigma_1 + \bm E_{11}) \bigr) + \bm Q^T \bigl( \bm E_{22} + \bm Q (\bm \Sigma_1 + \bm E_{11}) \bm P^T \bigr) \bm P \\
    &= (\bm Q^T \bm E_{22} - \bm Q^T \bm E_{21} \bm P^T) \bm P - \bm Q^T (\bm E_{22} \bm P + \bm Q \bm E_{12} \bm P) + \bm Q^T \bigl( \bm E_{22} + \bm Q (\bm \Sigma_1 + \bm E_{11}) \bm P^T \bigr) \bm P . \numberthis \label{equ:Pr153}
\end{align*}
Substituting (\ref{equ:Pr152}) and (\ref{equ:Pr153}) back into (\ref{equ:Pr151}), we obtain
\begin{align*}
    \P_r(\bm X + \bm \Delta) &= \bm X + \bm \Delta - \bm P_{\bm U_2} \bm \Delta \bm P_{\bm V_2} \\
    &\qquad + \bm U_1 \bm Q^T (\bm E_{21} \bm P^T - \bm E_{22}) \bm V_2^T + \bm U_2 (\bm E_{22} - \bm Q \bm E_{12}) \bm P \bm V_1^T + \bm U_2 (\bm E_{21} + \bm E_{22} \bm P + \bm Q \bm E_{12} \bm P) \bm P^T \bm V_2^T \\
    &\qquad + (\bm U_1 - \bm U_2 \bm Q) (\bm I_r + \bm Q^T \bm Q)^{-1} \cdot \Bigl( (\bm Q^T \bm E_{22} - \bm Q^T \bm E_{21} \bm P^T) \bm P - \bm Q^T (\bm E_{22} \bm P + \bm Q \bm E_{12} \bm P) \\
    &\qquad \qquad + \bm Q^T \bigl( \bm E_{22} + \bm Q (\bm \Sigma_1 + \bm E_{11}) \bm P^T \bigr) \bm P \Bigr) (\bm I_r + \bm P^T \bm P)^{-1} (\bm V_1^T + \bm P^T \bm V_2^T) . \numberthis \label{equ:Pr154}
\end{align*}
Since $\bm Q, \bm P, \bm E_{11}, \bm E_{12}, \bm E_{21}$, and $\bm E_{22}$ are first-order, and $(\bm I_r + \bm Q^T \bm Q)^{-1}$, $(\bm I_r + \bm P^T \bm P)^{-1}$ are zero-order in terms of $\norm{\bm \Delta}_F$, we can collect all the third-order terms on the RHS of (\ref{equ:Pr154}) and obtain
\begin{align} \label{equ:Pr155}
    \P_r(\bm X + \bm \Delta) = \bm X + \bm \Delta - \bm P_{\bm U_2} \bm \Delta \bm P_{\bm V_2} - \bm U_1 \bm Q^T \bm E_{22} \bm V_2^T + \bm U_2 \bm E_{22} \bm P \bm V_1^T + \bm U_2 \bm E_{21} \bm P^T \bm V_2^T + \bm \O(\norm{\bm \Delta}_F^3) .
\end{align}
Finally, the matrices $\bm Q$ and $\bm P$ in the second-order terms is eliminated by the following variant of (\ref{equ:QPS}):
\begin{align*}
    \bm Q &= - \bigl( \bm E_{21} + \bm Q \bm E_{21} \bm P - \bm E_{22} \bm P - \bm Q \bm E_{11} \bigr) \bm \Sigma_1^{-1} , \\
    \bm P^T &= \bm \Sigma_1^{-1} (\bm E_{12} + \bm Q^T \bm E_{21} \bm P^T - \bm Q^T \bm E_{22} - \bm E_{11} \bm P^T) .
\end{align*}
The substitution and collection of third-order terms on the RHS of (\ref{equ:Pr155}) yield
\begin{align*}
    \P_r(\bm X + \bm \Delta) &= \bm X + \bm \Delta - \bm P_{\bm U_2} \bm \Delta \bm P_{\bm V_2} + \bm U_1 \bm \Sigma_1^{-1} \bm E_{21}^T \bm E_{22} \bm V_2^T + \bm U_2 \bm E_{22} \bm E_{12}^T \bm \Sigma_1^{-1} \bm V_1^T + \bm U_2 \bm E_{21} \bm \Sigma_1^{-1} \bm E_{12} \bm V_2^T + \bm \O(\norm{\bm \Delta}_F^3) \\
    &= \bm X + \bm \Delta - \bm P_{\bm U_2} \bm \Delta \bm P_{\bm V_2} + \bm U_1 \bm \Sigma_1^{-1} \bm V_1^T \bm \Delta^T \bm U_2 \bm U_2^T \bm \Delta \bm V_2 \bm V_2^T \\
    &\qquad + \bm U_2 \bm U_2^T \bm \Delta \bm V_2 \bm V_2^T \bm \Delta^T \bm U_1 \bm \Sigma_1^{-1} \bm V_1^T + \bm U_2 \bm U_2^T \bm \Delta \bm V_1 \bm \Sigma_1^{-1} \bm U_1^T \bm \Delta \bm V_2 \bm V_2^T + \bm \O(\norm{\bm \Delta}_F^3) \\
    &= \bm X + \bm \Delta - \bm P_{\bm U_2} \bm \Delta \bm P_{\bm V_2} + \bm X^\dagger \bm \Delta^T P_{\bm U_2} \bm \Delta \bm P_{\bm V_2} + \bm P_{\bm U_2} \bm \Delta \bm P_{\bm V_2} \bm \Delta^T \bm X^\dagger + \bm P_{\bm U_2} \bm \Delta {(\bm X^\dagger)}^T \bm \Delta \bm P_{\bm V_2} + \bm \O(\norm{\bm \Delta}_F^3) .
\end{align*}
This completes our proof of the theorem.

\section{Proof of Lemma~\ref{lem:bound0}}
\label{apdx:bound0}

By the triangle inequality, we have
\begin{align*} 
    \norm{\P_r (\bm X + \bm \Delta) - (\bm X + \bm \Delta) + \bm P_{\bm U_2} \bm \Delta \bm P_{\bm V_2}}_F &\leq \norm{\P_r (\bm X + \bm \Delta) - (\bm X + \bm \Delta)}_F + \norm{\bm P_{\bm U_2} \bm \Delta \bm P_{\bm V_2}}_F . \numberthis \label{equ:Rm0}
\end{align*}
The first term on the RHS of (\ref{equ:Rm0}) can be bounded as follows. Since $\tilde{\bm X} = \bm X + \bm \Delta$, applying the norm absolute homogeneity property yields
\begin{align} \label{equ:Rm3}
    \norm{\P_r (\bm X + \bm \Delta) - (\bm X + \bm \Delta)}_F = \norm{\P_r (\tilde{\bm X}) - \tilde{\bm X}}_F = \norm{\tilde{\bm X} - \P_r (\tilde{\bm X})}_F .
\end{align}
From Lemmas~\ref{lem:Pr} and \ref{lem:semi}, we obtain
\begin{align} \label{equ:Rm4}
    \norm{\tilde{\bm X} - \P_r (\tilde{\bm X})}_F &= \norm{\tilde{\bm U}_{2} \tilde{\bm \Sigma}_2 \tilde{\bm V}_{2}^T}_F = \norm{\tilde{\bm \Sigma}_2}_F .
\end{align}
Since $\tilde{\bm \Sigma}_2$ is a submatrix of $\tilde{\bm \Sigma}$ containing $n-r$ small singular values of $\tilde{\bm X}$ in the diagonal, it holds that 
\begin{align} \label{equ:Rm5}
    \norm{\tilde{\bm \Sigma}_2}_F \leq \norm{\tilde{\bm \Sigma}}_F = \norm{\tilde{\bm X}}_F .
\end{align}
Additionally, using the triangle inequality we can bound $\norm{\tilde{\bm X}}_F$ by
\begin{align} \label{equ:Rm1}
    \norm{\tilde{\bm X}}_F = \norm{\bm X + \bm \Delta}_F \leq \norm{\bm X}_F + \norm{\bm \Delta}_F .
\end{align}
From (\ref{equ:Rm3}), (\ref{equ:Rm4}), (\ref{equ:Rm5}), and (\ref{equ:Rm1}), we have
\begin{align} \label{equ:Rm6}
    \norm{\P_r (\bm X + \bm \Delta) - (\bm X + \bm \Delta)}_F \leq \norm{\bm X}_F + \norm{\bm \Delta}_F .
\end{align}
On the other hand, it follows from Lemma~\ref{lem:semi} that the second term on the RHS of (\ref{equ:Rm0}) satisfies
\begin{align} \label{equ:Rm2}
    \norm{\bm P_{\bm U_2} \bm \Delta \bm P_{\bm V_2}}_F \leq \norm{\bm \Delta}_F .
\end{align}
Substituting inequalities (\ref{equ:Rm6}) and (\ref{equ:Rm2}) into (\ref{equ:Rm0}) completes the proof of the lemma.

\section{Proof of Theorem~\ref{THEO:RANKOP}}
\label{apdx:rankop}

\hlnew{The following proof is developed for the case of a rank-$r$ matrix $\bm X$. We first derive the proof of} (\ref{equ:boundR_1}) \hlnew{and then use this result to prove} (\ref{equ:boundR_2}).

\subsection{Proof of the bound in (\ref{equ:boundR_1})}

Our goal is to prove that the residual in (\ref{equ:rankop}) is always bounded by
\begin{align*}
    \norm{\bm R_{\bm X}(\bm \Delta)}_F \leq  \frac{c}{\sigma_r} \norm{\bm \Delta}_F^2 , \qquad \text{ for some } 1 + 1/\sqrt{2} \leq c \leq 4(1+\sqrt{2}) . 
\end{align*}
\hlnew{Let us begin with the upper bound on $c$ by showing that }
\begin{align} \label{equ:upper_c}
    \norm{\bm R_{\bm X}(\bm \Delta)}_F &\leq \frac{4(1+\sqrt{2})}{\sigma_r} \norm{\bm \Delta}_F^2 .
\end{align} 
Rearranging terms in (\ref{equ:rankop}) and replacing $\bm X + \bm \Delta$ by $\tilde{\bm X}$, we have 
\begin{align} \label{equ:R1}
    \bm R_{\bm X}(\bm \Delta) &= \P_r (\tilde{\bm X}) - \tilde{\bm X} + \bm P_{\bm U_2} \bm \Delta \bm P_{\bm V_2} .
\end{align}
Using the singular subspace decomposition in Definition~\ref{def:SVD} with descending order of singular values $\tilde{\sigma}_1 \geq \tilde{\sigma}_2 \ldots \geq \tilde{\sigma}_n$, let us decompose $\tilde{\bm X}$ as follows 
\begin{align} \label{equ:R11}
    \tilde{\bm X} = \tilde{\bm U}_1 \tilde{\bm \Sigma}_1 \tilde{\bm V}_1^T + \tilde{\bm U}_2 \tilde{\bm \Sigma}_2 \tilde{\bm V}_2^T .
\end{align}
Since in this theorem we consider perturbations of any magnitude, $\tilde{\bm X}$ can take any value including the case in which $\tilde{\sigma}_r=\tilde{\sigma}_{r+1}$ and the decomposition (\ref{equ:R11}) may not be unique. Nevertheless, the proof holds for any valid choice of singular subspace decomposition. From such a choice in (\ref{equ:R11}), $\P_r (\tilde{\bm X})$ is well-defined as: \hlnew{$\P_r (\tilde{\bm X}) = \tilde{\bm U}_1 \tilde{\bm \Sigma}_1 \tilde{\bm V}_1^T$.}
Substituting $\tilde{\bm X} = \bm X + \bm \Delta$ into (\ref{equ:G2}) and using the fact that $\bm P_{\bm U_2} \bm X = \bm 0$ and $\bm X \bm P_{\bm V_2} = \bm 0$, we obtain
\begin{align*}
    \bm R_{\bm X}(\bm \Delta) &= - \bm \delta_{\bm P_{\bm U_2}} \bm X \bm \delta_{\bm P_{\bm V_2}} - \bm P_{\bm U_2} \bm \Delta \bm \delta_{\bm P_{\bm V_2}} - \bm \delta_{\bm P_{\bm U_2}} \bm \Delta \bm P_{\bm V_2} - \bm \delta_{\bm P_{\bm U_2}} \bm \Delta \bm \delta_{\bm P_{\bm V_2}} \\
    &= - \bm \delta_{\bm P_{\bm U_2}} \bm X \bm \delta_{\bm P_{\bm V_2}} - \bm P_{\bm U_2} \bm \Delta \bm \delta_{\bm P_{\bm V_2}} - \bm \delta_{\bm P_{\bm U_2}} \bm \Delta \bm P_{\tilde{\bm V}_2} . \numberthis \label{equ:R13}
\end{align*}
Here, from Lemma~\ref{lem:trivial}, we can replace $\bm X = \bm X (\bm X^\dagger)^T \bm X$ in the first term on the RHS of (\ref{equ:R13}) and obtain
\begin{align} \label{equ:R14}
    \bm R_{\bm X}(\bm \Delta) = - (\bm \delta_{\bm P_{\bm U_2}} \bm X) (\bm X^\dagger)^T (\bm X \bm \delta_{\bm P_{\bm V_2}}) - \bm P_{\bm U_2} \bm \Delta \bm \delta_{\bm P_{\bm V_2}} - \bm \delta_{\bm P_{\bm U_2}} \bm \Delta \bm P_{\tilde{\bm V}_2} .
\end{align}
Taking the Frobenius norm and using its absolute homogeneity property, (\ref{equ:R14}) becomes
\begin{align*}
    \norm{\bm R_{\bm X}(\bm \Delta)}_F = \norm{(\bm \delta_{\bm P_{\bm U_2}} \bm X) (\bm X^\dagger)^T (\bm X \bm \delta_{\bm P_{\bm V_2}}) + \bm P_{\bm U_2} \bm \Delta \bm \delta_{\bm P_{\bm V_2}} + \bm \delta_{\bm P_{\bm U_2}} \bm \Delta \bm P_{\tilde{\bm V}_2}}_F .
\end{align*}
By the triangle inequality, the norm of $\bm R_{\bm X}(\bm \Delta)$ is then bounded by
\begin{align}\label{equ:R2}
    \norm{\bm R_{\bm X}(\bm \Delta)}_F &\leq \norm{(\bm \delta_{\bm P_{\bm U_2}} \bm X) (\bm X^\dagger)^T (\bm X \bm \delta_{\bm P_{\bm V_2}})}_F + \norm{\bm P_{\bm U_2} \bm \Delta \bm \delta_{\bm P_{\bm V_2}}}_F + \norm{\bm \delta_{\bm P_{\bm U_2}} \bm \Delta \bm P_{\tilde{\bm V}_2}}_F . 
\end{align}
Let us proceed to upper-bound $\norm{\bm R_{\bm X}(\bm \Delta)}_F$ by finding the upper bounds for each of the three terms on the RHS of (\ref{equ:R2}) with respect to $\norm{\bm \Delta}_F^2$. Our proof technique utilizes the following lemmas.

\begin{lemma} \label{lem:Mdp}
$\max \bigl\{ \norm{\bm \delta_{\bm P_{\bm U_2}} \bm X}_F, \norm{\bm X \bm \delta_{\bm P_{\bm V_2}}}_F \bigr\} \leq 2 \norm{\bm \Delta}_F.$
\end{lemma}

\begin{lemma} \label{lem:Pdp}
$\max \bigl\{ \norm{\bm P_{\bm U_2} \bm \delta_{\bm P_{\bm U_2}} \bm \Delta}_F, \norm{\bm \Delta \bm \delta_{\bm P_{\bm V_2}} \bm P_{\bm V_2}}_F \bigr\} \leq \frac{2}{\sigma_r} \norm{\bm \Delta}_F^2.$
\end{lemma}
\noindent The proofs of Lemmas~\ref{lem:Mdp} and \ref{lem:Pdp} are given at the end of this subsection. Let us proceed with the task of bounding the first term in (\ref{equ:R2}). Applying Lemma~\ref{lem:norm} twice and using the fact that $\norm {\bm X^{\dagger}}_2=1/\sigma_r$, we have
\begin{align*}
    \norm{(\bm \delta_{\bm P_{\bm U_2}} \bm X) (\bm X^\dagger)^T (\bm X \bm \delta_{\bm P_{\bm V_2}})}_F &\leq \frac{1}{\sigma_r} \norm{\bm \delta_{\bm P_{\bm U_2}} \bm X}_F \norm{\bm X \bm \delta_{\bm P_{\bm V_2}}}_F . \numberthis \label{equ:R1st0}
\end{align*}
By Lemma~\ref{lem:Mdp}, the terms $\norm{\bm \delta_{\bm P_{\bm U_2}} \bm X}_F$ and $\norm{\bm X \bm \delta_{\bm P_{\bm V_2}}}_F$ can each be bounded by $2 \norm{\bm \Delta}_F$. Applying the upper bounds on the RHS of  (\ref{equ:R1st0}), we obtain the  following bound on the first term in (\ref{equ:R2}):
\begin{align} \label{equ:R1st}
    \norm{(\bm \delta_{\bm P_{\bm U_2}} \bm X) (\bm X^\dagger)^T (\bm X \bm \delta_{\bm P_{\bm V_2}})}_F \le \frac{4}{\sigma_r} \norm{\bm \Delta}_F^2 .
\end{align}
Next, we shall bound the second term in (\ref{equ:R2}), i.e., $\norm{\bm P_{\bm U_2} \bm \Delta \bm \delta_{\bm P_{\bm V_2}}}_F$. From Lemma~\ref{lem:semi}, we have
\begin{align} \label{equ:R2nd1}
    \norm{\bm P_{\bm U_2} \bm \Delta \bm \delta_{\bm P_{\bm V_2}}}_F &\leq \norm{\bm \Delta \bm \delta_{\bm P_{\bm V_2}}}_F .
\end{align}
Since $\bm P_{\bm V_1} + \bm P_{\bm V_2} = \bm I_n$, the matrix on the RHS of (\ref{equ:R2nd1}) can be expanded as the sum of two orthogonal terms:
\begin{align*}
    \bm \Delta \bm \delta_{\bm P_{\bm V_2}} = \bm \Delta \bm \delta_{\bm P_{\bm V_2}} (\bm P_{\bm V_1} + \bm P_{\bm V_2}) = \bm \Delta \bm \delta_{\bm P_{\bm V_2}} \bm P_{\bm V_1} + \bm \Delta \bm \delta_{\bm P_{\bm V_2}} \bm P_{\bm V_2} .
\end{align*}
Notice that $\bm P_{\bm V_1}$ and $\bm P_{\bm V_2}$ are orthogonal. By Lemma~\ref{lem:Pythagoras}, we have 
\begin{align*}
    \norm{\bm \Delta \bm \delta_{\bm P_{\bm V_2}}}_F^2 &= \norm{\bm \Delta \bm \delta_{\bm P_{\bm V_2}} \bm P_{\bm V_1}}_F^2 + \norm{\bm \Delta \bm \delta_{\bm P_{\bm V_2}} \bm P_{\bm V_2}}_F^2 \\
    &= \norm{\bm \Delta \bm \delta_{\bm P_{\bm V_2}} \bm X^T \bm X^\dagger}_F^2 + \norm{\bm \Delta \bm \delta_{\bm P_{\bm V_2}} \bm P_{\bm V_2}}_F^2 && \text{(since } \bm P_{\bm V_1} = \bm X^T \bm X^\dagger \text{)} \\
    &= \norm{\bm \Delta (\bm X \bm \delta_{\bm P_{\bm V_2}})^T \bm X^\dagger}_F^2 + \norm{\bm \Delta \bm \delta_{\bm P_{\bm V_2}} \bm P_{\bm V_2}}_F^2 . \numberthis \label{equ:R2nd2}
\end{align*}
Each term on the RHS of (\ref{equ:R2nd2}) can be bounded as follows. Applying Lemma~\ref{lem:norm} twice, we initially bound the first term on the RHS of (\ref{equ:R2nd2}) as follows:
\begin{align*}
    \norm{\bm \Delta (\bm X \bm \delta_{\bm P_{\bm V_2}})^T \bm X^\dagger}_F &\leq \frac{1}{\sigma_r} \norm{\bm \Delta}_F \norm{\bm X \bm \delta_{\bm P_{\bm V_2}}}_F .
\end{align*}
By Lemma~\ref{lem:Mdp}, we upper-bound $\norm{\bm X \bm \delta_{\bm P_{\bm V_2}}}_F$ by $2\norm{\bm \Delta}_F$ and obtain the bound on the first term on the RHS of (\ref{equ:R2nd2}):
\begin{align} \label{equ:R2nd3}
    \norm{\bm \Delta (\bm X \bm \delta_{\bm P_{\bm V_2}})^T \bm X^\dagger}_F \leq \frac{2}{\sigma_r} \norm{\bm \Delta}_F^2 .
\end{align}
To bound the second term on the RHS of (\ref{equ:R2nd2}), we apply  Lemma~\ref{lem:Pdp} and obtain 
\begin{align} \label{equ:R2nd4}
    \norm{\bm \Delta \bm \delta_{\bm P_{\bm V_2}} \bm P_{\bm V_2}}_F^2 \leq \frac{4}{\sigma_r^2} \norm{\bm \Delta}_F^4 .
\end{align}
\hlnew{Substituting the bounds from} (\ref{equ:R2nd3}) and (\ref{equ:R2nd4}) back into the RHS of (\ref{equ:R2nd2}), we have
\begin{align*}
    \norm{\bm \Delta \bm \delta_{\bm P_{\bm V_2}}}_F^2 &\leq \Bigl( \frac{2}{\sigma_r} \norm{\bm \Delta}_F^2 \Bigr)^2 + \frac{4}{\sigma_r^2} \norm{\bm \Delta}_F^4 = \frac{8}{\sigma_r^2} \norm{\bm \Delta}_F^4.
\end{align*}
Taking the square root of the last result and substituting it back to (\ref{equ:R2nd1})  yields
\begin{align} \label{equ:R2nd}
    \norm{\bm P_{\bm U_2} \bm \Delta \bm \delta_{\bm P_{\bm V_2}}}_F \leq \frac{2\sqrt{2}}{\sigma_r} \norm{\bm \Delta}_F^2 .
\end{align}
This offers a bound on the second term on the RHS of (\ref{equ:R2}).
Similarly, we bound the third term on the RHS of (\ref{equ:R2}) by 
\begin{align} \label{equ:R3rd}
    \norm{\bm \delta_{\bm P_{\bm U_2}} \bm \Delta \bm P_{\tilde{\bm V}_2}}_F \leq \frac{2\sqrt{2}}{\sigma_r} \norm{\bm \Delta}_F^2 .
\end{align}
Finally, summing up (\ref{equ:R1st}), (\ref{equ:R2nd}), and (\ref{equ:R3rd}), and substituting back into (\ref{equ:R2}), we obtain (\ref{equ:upper_c}) \hlnew{and thereby completes the first part of the proof.}

\hlnew{For the second part of the proof, we show that $c \geq 1 + 1/\sqrt{2}$ by constructing a perturbation $\bm \Delta$ such that the ratio $\norm{\bm R_{\bm X}(\bm \Delta)}_F/\norm{\bm \Delta}_F^2$ approaches $(1 + 1/\sqrt{2})/\sigma_r$.}
Consider perturbations of form
\begin{align} \label{equ:c1_0}
    \bm \Delta = (\sigma - \sigma_r - \epsilon) \bm u_r \bm v_r^T + \sigma \bm u_{r+1} \bm v_{r+1}^T, \qquad \text{for } 0<\epsilon<\sigma<\sigma_r .
\end{align}
Since $\bm u_r \bm v_r^T$ and $\bm u_{r+1} \bm v_{r+1}^T$ are orthogonal, we can compute the norm of $\bm \Delta$ using Lemma~\ref{lem:Pythagoras}:
\begin{align*} 
    \norm{\bm \Delta}_F^2 &= (\sigma - \sigma_r - \epsilon)^2 \norm{\bm u_r \bm v_r^T}_F^2 + \sigma^2 \norm{\bm u_{r+1} \bm v_{r+1}^T}_F^2 \\
    &= (\sigma - \sigma_r - \epsilon)^2 + \sigma^2 , \numberthis \label{equ:c1_1} 
\end{align*}
where the second equality uses $\bm u_r \bm v_r^T = \bm u_r \otimes \bm v_r^T$ and Lemma~\ref{lem:kron}-3.
Using the SVD of ${\bm X}$ and the definition of ${\bm \Delta}$ in (\ref{equ:c1_0}), we have
\begin{align*} 
    \bm X + \bm \Delta &= \sum_{i=1}^{r} \sigma_i \bm u_i \bm v_i^T + (\sigma - \sigma_r - \epsilon) \bm u_r \bm v_r^T + \sigma \bm u_{r+1} \bm v_{r+1}^T \\
    &= \sum_{i=1}^{r-1} \sigma_i \bm u_i \bm v_i^T + (\sigma-\epsilon) \bm u_r \bm v_r^T + \sigma \bm u_{r+1} \bm v_{r+1}^T . \numberthis \label{equ:c1_2}
\end{align*}
After perturbation, the $r$-th singular value of ${\bm X}$ is changed from $\sigma_r$ to $\sigma-\epsilon$ and the $r+1$-th changes from $0$ to $\sigma$, thereby making the singular value corresponding to $\bm u_{r+1} \bm v_{r+1}^T$ larger than the singular value associated with $\bm u_r \bm v_r^T$. Thus, the $r$-TSVD of $\bm X + \bm \Delta$ is given by
\begin{align} \label{equ:c1_3}
     \P_r (\bm X + \bm \Delta) &= \sum_{i=1}^{r-1} \sigma_i \bm u_i \bm v_i^T + \sigma \bm u_{r+1} \bm v_{r+1}^T .
\end{align}
On the other hand, since $\bm P_{\bm U_2} = \sum_{i=r+1}^m \bm u_i \bm u_i^T$ and $\bm P_{\bm V_2} = \sum_{i=r+1}^n \bm v_i \bm v_i^T$, we have
\begin{align*} 
    \bm P_{\bm U_2} \bm \Delta \bm P_{\bm V_2} &= \Bigl(\sum_{i=r+1}^m \bm u_i \bm u_i^T \Bigr) \Bigl( (\sigma - \sigma_r - \epsilon) \bm u_r \bm v_r^T + \sigma \bm u_{r+1} \bm v_{r+1}^T \Bigr) \Bigl(\sum_{i=r+1}^n \bm v_i \bm v_i^T \Bigr) = \sigma \bm u_{r+1} \bm v_{r+1}^T , \numberthis \label{equ:c1_31}
\end{align*}
where the second equality stems from the fact that
\begin{align*}
         \bm u_i^T \bm u_j = \bm v_i^T \bm v_j = \begin{cases} 1 & \text{if } i=j, \\ 0 & \text{if } i \neq j . \end{cases} 
\end{align*}
Substituting (\ref{equ:c1_2}), (\ref{equ:c1_3}), and (\ref{equ:c1_31}) into (\ref{equ:R1}), we obtain 
\begin{align*}
    \bm R_{\bm X}(\bm \Delta) &= \biggl( \sum_{i=1}^{r-1} \sigma_i \bm u_i \bm v_i^T + \sigma \bm u_{r+1} \bm v_{r+1}^T \biggr) - \biggl( \sum_{i=1}^{r-1} \sigma_i \bm u_i \bm v_i^T + (\sigma-\epsilon) \bm u_r \bm v_r^T + \sigma \bm u_{r+1} \bm v_{r+1}^T \biggr) + \sigma \bm u_{r+1} \bm v_{r+1}^T \\
    &= -(\sigma-\epsilon) \bm u_r \bm v_r^T + \sigma \bm u_{r+1} \bm v_{r+1}^T .
\end{align*}
Similar to (\ref{equ:c1_1}), one can compute the norm of the residual by
\begin{align} \label{equ:c1_4}
    \norm{\bm R_{\bm X}(\bm \Delta)}_F = \sqrt{(\sigma-\epsilon)^2 + \sigma^2} .
\end{align}
From (\ref{equ:c1_1}) and (\ref{equ:c1_4}), we have
\begin{align*}
    \frac{\norm{\bm R_{\bm X}(\bm \Delta)}_F}{\norm{\bm \Delta}_F^2} = \frac{\sqrt{(\sigma-\epsilon)^2 + \sigma^2}}{(\sigma_r + \epsilon - \sigma)^2 + \sigma^2} .
\end{align*}
Now maximizing over $\sigma$ while taking $\epsilon$ to $0$ gives us a lower bound on $c$:
\begin{align*}
    \frac{c}{\sigma_r} &= \sup_{\bm \Delta \in \R^{m \times n}} \frac{\norm{\bm R_{\bm X}(\bm \Delta)}_F}{\norm{\bm \Delta}_F^2} \\
    &\geq \max_{0<\sigma<\sigma_r} ~\lim_{\epsilon \to 0^+} \frac{\sqrt{(\sigma-\epsilon)^2 + \sigma^2}}{(\sigma_r + \epsilon - \sigma)^2 + \sigma^2} \\
    &= \max_{0<\sigma<\sigma_r} \frac{\sigma \sqrt{2}}{(\sigma_r - \sigma)^2 + \sigma^2} . \numberthis \label{equ:c1s_1}
\end{align*}
The maximization can be obtained at $\sigma = \sigma_r/\sqrt{2}$. Therefore, substituting back into (\ref{equ:c1s_1}) yields $c \geq 1 + 1/\sqrt{2}$.
This completes our proof of the first half of Theorem~\ref{THEO:RANKOP}. We recall from Remark~\ref{rem:c} \hlnew{that we conjecture the structure of $\bm \Delta$ given in} (\ref{equ:c1_0}) \hlnew{yields the maximizer of $\norm{\bm R_{\bm X}(\bm \Delta)}_F / \norm{\bm \Delta}_F^2$.}

\subsubsection{Proof of Lemma~\ref{lem:Mdp}}
Let us rewrite $\bm \delta_{\bm P_{\bm U_2}} \bm X = \bm P_{\tilde{\bm U}_2} \bm X - \bm P_{\bm U_2} \bm X$.
Since $\bm P_{\bm U_2} \bm X = \bm 0$, we obtain
\begin{align*}
    \bm \delta_{\bm P_{\bm U_2}} \bm X &= \bm P_{\tilde{\bm U}_2} \bm X \numberthis \label{equ:MdP0} \\
    &= \bm P_{\tilde{\bm U}_2} (\tilde{\bm X} - \bm \Delta) && \text{(since } \tilde{\bm X} = \bm X + \bm \Delta \text{)} \\
    &= \tilde{\bm U}_2 \tilde{\bm U}_2^T \tilde{\bm X} - \bm P_{\tilde{\bm U}_2} \bm \Delta .
\end{align*}
Substituting $\tilde{\bm X} = \tilde{\bm U}_1 \tilde{\bm \Sigma}_1 \tilde{\bm V}_1^T + \tilde{\bm U}_2 \tilde{\bm \Sigma}_2 \tilde{\bm V}_2^T$ yields
\begin{align*}
    \bm \delta_{\bm P_{\bm U_2}} \bm X &= \tilde{\bm U}_2 \tilde{\bm U}_2^T \bigl(\tilde{\bm U}_1 \tilde{\bm \Sigma}_1 \tilde{\bm V}_1^T + \tilde{\bm U}_2 \tilde{\bm \Sigma}_2 \tilde{\bm V}_2^T\bigr) - \bm P_{\tilde{\bm U}_2} \bm \Delta \\
    &= \tilde{\bm U}_2 \tilde{\bm \Sigma}_2 \tilde{\bm V}_2^T - \bm P_{\tilde{\bm U}_2} \bm \Delta , 
\end{align*}
where in the last equality we use the fact that $\tilde{\bm U}_2^T \tilde{\bm U}_1 = \bm 0$ and $\tilde{\bm U}_2^T \tilde{\bm U}_2 = \bm I_m$.
Therefore,
\begin{align} \label{equ:MdP1}
    \norm{\bm \delta_{\bm P_{\bm U_2}} \bm X}_F = \norm{\tilde{\bm U}_2 \tilde{\bm \Sigma}_2 \tilde{\bm V}_2^T - \bm P_{\tilde{\bm U}_2} \bm \Delta}_F .
\end{align}
By the triangle inequality and the absolute homogeneity, (\ref{equ:MdP1}) implies
\begin{align} \label{equ:MdP2}
    \norm{\bm \delta_{\bm P_{\bm U_2}} \bm X}_F &\leq \norm{\tilde{\bm U}_2 \tilde{\bm \Sigma}_2 \tilde{\bm V}_2^T}_F + \norm{\bm P_{\tilde{\bm U}_2} \bm \Delta}_F .
\end{align}
We shall bound each term on the RHS of (\ref{equ:MdP2}) as follows. First, using Lemma~\ref{lem:semi}, we can  remove the semi-orthogonal matrices from within the Frobenius norm without changing the value of the norm:
\begin{align*}
    \norm{\tilde{\bm U}_2 \tilde{\bm \Sigma}_2 \tilde{\bm V}_2^T}_F = \norm{\tilde{\bm \Sigma}_2 \tilde{\bm V}_2^T}_F = \norm{\tilde{\bm \Sigma}_2}_F .
\end{align*}
Since $\bm \Sigma_2 = \bm 0$, we further obtain
\begin{align} \label{equ:MdP3}
    \norm{\tilde{\bm U}_2 \tilde{\bm \Sigma}_2 \tilde{\bm V}_2^T}_F = \norm{\tilde{\bm \Sigma}_2 - \bm \Sigma_2}_F .
\end{align}
In addition, recall that $\tilde{\bm \Sigma}_2$ and $\bm \Sigma_2$ are sub-matrices of $\tilde{\bm \Sigma}$ and $\bm \Sigma$, respectively. Thus, 
\begin{align} \label{equ:MdP40}
    \norm{\tilde{\bm \Sigma}_2 - \bm \Sigma_2}_F &\leq \norm{\tilde{\bm \Sigma} - \bm \Sigma}_F .
\end{align}
Moreover, by Mirsky's inequality in Proposition~\ref{pro:s}, we have
\begin{align*}
    \norm{\tilde{\bm \Sigma} - \bm \Sigma}_F = \sqrt{\sum_{i=1}^n (\tilde{\sigma}_i - \sigma_i)^2} &\leq \norm{\bm \Delta}_F . \numberthis \label{equ:MdP4}
\end{align*}
From (\ref{equ:MdP3}), (\ref{equ:MdP40}), and (\ref{equ:MdP4}), it follows that
\begin{align} \label{equ:MdP5a}
    \norm{\tilde{\bm U}_2 \tilde{\bm \Sigma}_2 \tilde{\bm V}_2^T}_F \leq \norm{\bm \Delta}_F .
\end{align}
Next, the second term on the RHS of (\ref{equ:MdP2}), by Lemma~\ref{lem:semi}, is bounded by
\begin{align} \label{equ:MdP5}
    \norm{\bm P_{\tilde{\bm U}_2} \bm \Delta}_F \leq \norm{\bm \Delta}_F .
\end{align}
Summing up (\ref{equ:MdP5a}) and (\ref{equ:MdP5}), and combining the resulting inequality with (\ref{equ:MdP2}), we conclude that
\begin{align*} 
    \norm{\bm \delta_{\bm P_{\bm U_2}} \bm X}_F &\leq 2 \norm{\bm \Delta}_F . 
\end{align*}
The proof of $\norm{\bm X \bm \delta_{\bm P_{\bm V_2}}}_F \leq 2 \norm{\bm \Delta}_F$ follows a similar derivation.

\subsubsection{Proof of Lemma~\ref{lem:Pdp}}

In this subsection, we shall show that $\norm{\bm P_{\bm U_2} \bm \delta_{\bm P_{\bm U_2}} \bm \Delta}_F \leq \frac{2}{\sigma_r} \norm{\bm \Delta}_F^2$. The proof of $\norm{\bm \Delta \bm \delta_{\bm P_{\bm V_2}} \bm P_{\bm V_2}}_F \leq \frac{2}{\sigma_r} \norm{\bm \Delta}_F^2$ can be derived similarly. 
Since Definition~\ref{def:Puv} implies $\bm \delta_{\bm P_{\bm U_2}} = \bm P_{\tilde{\bm U}_2} - \bm P_{\bm U_2} = \bm P_{\bm U_1} - \bm P_{\tilde{\bm U}_1}$, we have
\begin{align*} 
    \bm P_{\bm U_2} \bm \delta_{\bm P_{\bm U_2}} \bm \Delta &= \bm P_{\bm U_2} (\bm P_{\bm U_1} - \bm P_{\tilde{\bm U}_1}) \bm \Delta \\
    &= -\bm P_{\bm U_2} \bm P_{\tilde{\bm U}_1} \bm \Delta , \numberthis \label{equ:Pdp_1}
\end{align*}
where the second equality is due to $\bm P_{\bm U_2} \bm P_{\bm U_1} = \bm 0$ (see Lemma~\ref{lem:trivial}).
It is now sufficient to bound the norm of $\bm P_{\bm U_2} \bm P_{\tilde{\bm U}_1} \bm \Delta$ by $\frac{2}{\sigma_r} \norm{\bm \Delta}_F^2$. Let us consider two cases:
\begin{itemize}
    \item If $\norm{\bm \Delta}_2 \geq \sigma_r/2$, then applying Lemma~\ref{lem:semi}-2 twice yields
    \begin{align*}
        \norm{\bm P_{\bm U_2} \bm P_{\tilde{\bm U}_1} \bm \Delta}_F &\leq \norm{\bm \Delta}_F . \numberthis \label{equ:Pdp_2}
    \end{align*}
    Since $\norm{\bm \Delta}_F \geq \norm{\bm \Delta}_2 \geq \sigma_r/2$, multiplying both sides by $\frac{2}{\sigma_r} \norm{\bm \Delta}_F$ yields
    \begin{align*}
        \norm{\bm \Delta}_F &\leq \frac{2}{\sigma_r} \norm{\bm \Delta}_F^2 . \numberthis \label{equ:Pdp_3}
    \end{align*}
    From (\ref{equ:Pdp_2}) and (\ref{equ:Pdp_3}), we obtain $\norm{\bm P_{\bm U_2} \bm P_{\tilde{\bm U}_1} \bm \Delta}_F \leq \frac{2}{\sigma_r} \norm{\bm \Delta}_F^2$.
    
    \item If $\norm{\bm \Delta}_2 < \sigma_r/2$, we need to use a different approach as follows. First, from Lemma~\ref{lem:norm}, we have
    \begin{align} \label{equ:Pdp_9}
        \norm{\bm P_{\bm U_2} \bm P_{\tilde{\bm U}_1} \bm \Delta}_F \leq \norm{\bm P_{\bm U_2} \bm P_{\tilde{\bm U}_1}}_2 \norm{\bm \Delta}_F .
    \end{align}
    Let us examine the product $\bm P_{\bm U_2} \bm P_{\tilde{\bm U}_1}$.
    Let $\tilde{\bm X}_1 = \tilde{\bm U}_1 \tilde{\bm \Sigma}_1 \tilde{\bm V}_1^T$ and $\tilde{\bm X}_2 = \tilde{\bm X} - \tilde{\bm X}_1$. From Weyl's inequality \cite{weyl1912asymptotische}, we have
    \begin{align*}
        \abs{\tilde{\sigma}_i - \sigma_i} \leq \norm{\bm \Delta}_2 < \frac{\sigma_r}{2} \qquad \text{for } i=1,\ldots,n.
    \end{align*}
    Thus, for any $1 \leq i \leq r$, it holds that
    \begin{align} \label{equ:Pdp_4}
        \tilde{\sigma}_i > \sigma_i - \frac{\sigma_r}{2} \geq \sigma_r - \frac{\sigma_r}{2} = \frac{\sigma_r}{2} > 0 .
    \end{align}
    Therefore, $\tilde{\bm \Sigma}_1 = \diag(\tilde{\sigma}_1,\ldots,\tilde{\sigma}_r)$ is invertible. We can now denote the pseudo inverse of $\tilde{\bm X}_1$ by $\tilde{\bm X}_1^\dagger = \tilde{\bm U}_1 \tilde{\bm \Sigma}_1^{-1} \tilde{\bm V}_1^T$. We have 
    \begin{align*}
        \bm P_{\bm U_2} \bm P_{\tilde{\bm U}_1} &= \bm P_{\bm U_2} \tilde{\bm X}_1 (\tilde{\bm X}_1^\dagger)^T && \text{(since } \bm P_{\tilde{\bm U}_1} = \tilde{\bm X}_1 (\tilde{\bm X}_1^\dagger)^T \text{)} \\
        &= \bm P_{\bm U_2} (\tilde{\bm X} - \tilde{\bm X}_2) (\tilde{\bm X}_1^\dagger)^T \\
        &= \bm P_{\bm U_2} \tilde{\bm X} (\tilde{\bm X}_1^\dagger)^T && \text{(since } \tilde{\bm X}_2 (\tilde{\bm X}_1^\dagger)^T = \bm 0 \text{)} \\
        &= \bm P_{\bm U_2} (\bm X + \bm \Delta) (\tilde{\bm X}_1^\dagger)^T \\
        &= \bm P_{\bm U_2} \bm \Delta (\tilde{\bm X}_1^\dagger)^T . && \text{(since } \bm P_{\bm U_2} \bm X = \bm 0 \text{)} \numberthis \label{equ:Pdp_5}
    \end{align*}
    On the other hand, applying Lemmas~\ref{lem:semi} and \ref{lem:norm}, and the fact that $\norm{\bm X^\dagger}_2=1/\sigma_r$, we obtain
    \begin{align*}
        \norm{\bm P_{\bm U_2} \bm \Delta (\tilde{\bm X}_1^\dagger)^T}_F \leq \frac{1}{\tilde{\sigma}_r} \norm{\bm \Delta}_F . \numberthis \label{equ:Pdp_6}
    \end{align*}
    From (\ref{equ:Pdp_4}), we can bound  $\tilde{\sigma}_r$ by:
    \begin{align} \label{equ:Pdp_7}
        \tilde{\sigma}_r > \sigma_r - \frac{\sigma_r}{2} = \frac{\sigma_r}{2} .
    \end{align}
    From (\ref{equ:Pdp_5}), (\ref{equ:Pdp_6}), and (\ref{equ:Pdp_7}), we obtain
    \begin{align} \label{equ:Pdp_8}
        \norm{\bm P_{\bm U_2} \bm P_{\tilde{\bm U}_1}}_F = \norm{\bm P_{\bm U_2} \bm \Delta (\tilde{\bm X}_1^\dagger)^T}_F \leq \frac{1}{\tilde{\sigma}_r} \norm{\bm \Delta}_F < \frac{2}{\sigma_r} \norm{\bm \Delta}_F .
    \end{align}
    Finally, substituting (\ref{equ:Pdp_8}) back into (\ref{equ:Pdp_9}) immediately yields $\norm{\bm P_{\bm U_2} \bm P_{\tilde{\bm U}_1} \bm \Delta}_F < \frac{2}{\sigma_r} \norm{\bm \Delta}_F^2$.
\end{itemize}
Since in both cases $\norm{\bm P_{\bm U_2} \bm P_{\tilde{\bm U}_1} \bm \Delta}_F \leq \frac{2}{\sigma_r} \norm{\bm \Delta}_F^2$, we conclude from (\ref{equ:Pdp_1}) that $\norm{\bm P_{\bm U_2} \bm \delta_{\bm P_{\bm U_2}} \bm \Delta}_F \leq \frac{2}{\sigma_r} \norm{\bm \Delta}_F^2$ for any ${\bm \Delta}$.

\subsection{Proof of the bound in (\ref{equ:boundR_2})}
\label{apdx:R}

Taking Frobenius norm on both sides of equation (\ref{equ:R13}) and using its absolute homogeneity property, we obtain:
\begin{align}\label{equ:Rg11}
    \| \bm R_{\bm X}(\bm \Delta)\| = \| \bm \delta_{\bm P_{\bm U_2}} \bm X \bm \delta_{\bm P_{\bm V_2}} + \bm P_{\bm U_2} \bm \Delta \bm \delta_{\bm P_{\bm V_2}} + \bm \delta_{\bm P_{\bm U_2}} \bm \Delta \bm P_{\tilde{\bm V}_2} \|.
\end{align}
Applying the triangle inequality to the RHS of (\ref{equ:Rg11}), we have
\begin{align}\label{equ:Rg12}
    \norm{\bm R_{\bm X}(\bm \Delta)}_F &\leq \norm{\bm \delta_{\bm P_{\bm U_2}} \bm X \bm \delta_{\bm P_{\bm V_2}}}_F + \norm{\bm P_{\bm U_2} \bm \Delta \bm \delta_{\bm P_{\bm V_2}}}_F + \norm{\bm \delta_{\bm P_{\bm U_2}} \bm \Delta \bm P_{\tilde{\bm V}_2}}_F .
\end{align}
To bound the RHS of (\ref{equ:Rg12}), we proceed by bounding each of the terms on the RHS. The first term on the RHS of (\ref{equ:Rg12}) can be bounded as follows. From (\ref{equ:MdP0}), we have $\bm \delta_{\bm P_{\bm U_2}} \bm X \bm \delta_{\bm P_{\bm V_2}} = \bm P_{\tilde{ \bm U}_2} \bm X \bm \delta_{\bm P_{\bm V_2}}$. Using Lemmas~\ref{lem:semi} and \ref{lem:Mdp}, it follows that
\begin{align*}
    \norm{\bm \delta_{\bm P_{\bm U_2}} \bm X \bm \delta_{\bm P_{\bm V_2}}}_F &= \norm{\bm P_{\tilde{ \bm U}_2} \bm X \bm \delta_{\bm P_{\bm V_2}}}_F \\
    &\leq \norm{\bm X \bm \delta_{\bm P_{\bm V_2}}}_F \\
    &\leq 2 \norm{\bm \Delta}_F . \numberthis \label{equ:Rg13}
\end{align*}
Next, the second term on the RHS of (\ref{equ:Rg12}) can be rewritten as the sum of two orthogonal components
\begin{align*}
    \bm P_{\bm U_2} \bm \Delta \bm \delta_{\bm P_{\bm V_2}} = \bm P_{\bm U_2} \bm \Delta \bm \delta_{\bm P_{\bm V_2}} \bm P_{\bm V_1} + \bm P_{\bm U_2} \bm \Delta \bm \delta_{\bm P_{\bm V_2}} \bm P_{\bm V_2} .
\end{align*}
By Lemma~\ref{lem:Pythagoras}, we have
\begin{align} \label{equ:Rg14}
   \norm{\bm P_{\bm U_2} \bm \Delta \bm \delta_{\bm P_{\bm V_2}}}_F = \sqrt{\norm{\bm P_{\bm U_2} \bm \Delta \bm \delta_{\bm P_{\bm V_2}} \bm P_{\bm V_1}}_F^2 + \norm{\bm P_{\bm U_2} \bm \Delta \bm \delta_{\bm P_{\bm V_2}} \bm P_{\bm V_2}}_F^2} .
\end{align}
On the one hand, we consider the first term on the RHS of (\ref{equ:Rg14}). Since 
\begin{align*}
    \bm \delta_{\bm P_{\bm V_2}} \bm P_{\bm V_1} &= (\bm P_{\tilde{\bm V}_2} - \bm P_{\bm V_2}) \bm P_{\bm V_1} \\
    &= \bm P_{\tilde{\bm V}_2} \bm P_{\bm V_1} , && \text{(by Lemma~\ref{lem:trivial})}
\end{align*}
we obtain
\begin{align} \label{equ:Rg15}
    \norm{\bm P_{\bm U_2} \bm \Delta \bm \delta_{\bm P_{\bm V_2}} \bm P_{\bm V_1}}_F &= \norm{\bm P_{\bm U_2} \bm \Delta \bm P_{\tilde{\bm V}_2} \bm P_{\bm V_1}}_F .
\end{align}
Applying Lemma~\ref{lem:semi} to the RHS of (\ref{equ:Rg15}) in order to eliminate the three projection matrices, we obtain
\begin{align} \label{equ:Rg16}
    \norm{\bm P_{\bm U_2} \bm \Delta \bm \delta_{\bm P_{\bm V_2}} \bm P_{\bm V_1}}_F \leq \norm{\bm \Delta}_F .
\end{align}
Similarly, we have
\begin{align} \label{equ:Rg18}
    \norm{\bm P_{\bm U_2} \bm \Delta \bm \delta_{\bm P_{\bm V_2}} \bm P_{\bm V_2}}_F \leq \norm{\bm \Delta}_F .
\end{align}
Substituting (\ref{equ:Rg16}), and (\ref{equ:Rg18}) back into (\ref{equ:Rg14}), we have
\begin{align} \label{equ:Rg19}
    \norm{\bm P_{\bm U_2} \bm \Delta \bm \delta_{\bm P_{\bm V_2}}}_F  &\leq \sqrt{2} \norm{\bm \Delta}_F .
\end{align}
Similarly, we also obtain
\begin{align} \label{equ:Rg20}
   \norm{ \bm \delta_{\bm P_{\bm U_2}} \bm \Delta \bm P_{\tilde{\bm V}_2} }_F  &\leq \sqrt{2} \norm{\bm \Delta}_F .
\end{align}
Substituting (\ref{equ:Rg13}), (\ref{equ:Rg19}), and (\ref{equ:Rg20}) back into (\ref{equ:Rg12}), we obtain
\begin{align} \label{equ:Rg21}
\norm{\bm R_{\bm X}(\bm \Delta)}_F \leq 2(1+\sqrt{2}) \norm{\bm \Delta}_F .
\end{align}
The proof of (\ref{equ:boundR_2}) is concluded by taking the minimum between the bounds in (\ref{equ:Rg21}) and (\ref{equ:boundR_1}).

\section{Proof of Theorem \ref{THEO:RSUP}}
\label{apdx:RSUP}

Let us denote
\begin{align*}
    \bm R_{2\bm X}(\bm \Delta) = \bm X^\dagger \bm \Delta^T \bm P_{\bm U_2} \bm \Delta \bm P_{\bm V_2} + \bm P_{\bm U_2} \bm \Delta \bm P_{\bm V_2} \bm \Delta^T \bm X^\dagger + \bm P_{\bm U_2} \bm \Delta {(\bm X^\dagger)}^T \bm \Delta \bm P_{\bm V_2} .
\end{align*}
It is straightforward to verify from (\ref{equ:rankop2}) that $\bm R_{\bm X}(\bm \Delta) = \bm R_{2\bm X}(\bm \Delta) + \bm \O(\norm{\bm \Delta}_F^3)$. Thus, 
\begin{align} \label{equ:Rg299}
    \lim_{\epsilon \to 0^+} \sup_{\norm{\bm \Delta}_F = \epsilon} \frac{\norm{\bm R_{\bm X}(\bm \Delta) - \bm R_{2\bm X}(\bm \Delta)}_F}{\norm{\bm \Delta}_F^2} = 0 .
\end{align}

\begin{lemma} \label{lem:max}
Let $f$ and $g$ be some bounded real-valued functions defined on the set $\C$. 
Then it holds that
\begin{align*}
    \abs{\sup_{\bm x \in \C} f(\bm x) - \sup_{\bm x \in \C} g(\bm x)} \leq \sup_{\bm x \in \C} \abs{f(\bm x) - g(\bm x)} .
\end{align*}
\end{lemma}
\noindent Applying Lemma~\ref{lem:max} to (\ref{equ:Rg299}), we obtain
\begin{align} \label{equ:Rg22a}
    \abs{\sup_{\norm{\bm \Delta}_F = \epsilon} \frac{\norm{\bm R_{\bm X}(\bm \Delta)}_F}{\norm{\bm \Delta}_F^2} - \sup_{\norm{\bm \Delta}_F = \epsilon} \frac{\norm{\bm R_{2\bm X}(\bm \Delta)}_F}{\norm{\bm \Delta}_F^2}} \leq \sup_{\norm{\bm \Delta}_F = \epsilon} \abs{\frac{\norm{\bm R_{\bm X}(\bm \Delta)}_F - \norm{\bm R_{2\bm X}(\bm \Delta)}_F}{\norm{\bm \Delta}_F^2}} .
\end{align}
On the other hand, by the triangle inequality, we have
\begin{align} \label{equ:Rg22b}
    \abs{\norm{\bm R_{\bm X}(\bm \Delta)}_F - \norm{\bm R_{2\bm X}(\bm \Delta)}_F}\leq \norm{\bm R_{\bm X}(\bm \Delta) - \bm R_{2\bm X}(\bm \Delta)}_F .
\end{align}
From (\ref{equ:Rg22a}) and (\ref{equ:Rg22b}), it holds that
\begin{align} \label{equ:Rg22}
    \abs{\sup_{\norm{\bm \Delta}_F = \epsilon} \frac{\norm{\bm R_{\bm X}(\bm \Delta)}_F}{\norm{\bm \Delta}_F^2} - \sup_{\norm{\bm \Delta}_F = \epsilon} \frac{\norm{\bm R_{2\bm X}(\bm \Delta)}_F}{\norm{\bm \Delta}_F^2}} \leq \sup_{\norm{\bm \Delta}_F = \epsilon} \abs{\frac{\norm{\bm R_{\bm X}(\bm \Delta) - \bm R_{2\bm X}(\bm \Delta)}_F}{\norm{\bm \Delta}_F^2}} .
\end{align}
Thus, taking both sides of (\ref{equ:Rg22}) to the limit $\epsilon \to 0$ and rearranging terms yield
\begin{align*}
    \lim_{\epsilon \to 0^+} \sup_{\norm{\bm \Delta}_F = \epsilon} \frac{\norm{\bm R_{\bm X}(\bm \Delta)}_F}{\norm{\bm \Delta}_F^2} = \lim_{\epsilon \to 0^+} \sup_{\norm{\bm \Delta}_F = \epsilon} \frac{\norm{\bm R_{2\bm X}(\bm \Delta)}_F}{\norm{\bm \Delta}_F^2} .
\end{align*}
It now is sufficient to show that 
\begin{align} \label{equ:Rg25}
    \lim_{\epsilon \to 0^+} \sup_{\norm{\bm \Delta}_F = \epsilon} \frac{\norm{\bm R_{2\bm X}(\bm \Delta)}_F}{\norm{\bm \Delta}_F^2} = \frac{1}{\sigma_r\sqrt{3}} .
\end{align}
Indeed, due to the orthogonality among the addends, we have
\begin{align*}
    \norm{\bm R_{2\bm X}(\bm \Delta)}_F^2 &= \norm{\bm X^\dagger \bm \Delta^T \bm P_{\bm U_2} \bm \Delta \bm P_{\bm V_2} + \bm P_{\bm U_2} \bm \Delta \bm P_{\bm V_2} \bm \Delta^T \bm X^\dagger + \bm P_{\bm U_2} \bm \Delta {(\bm X^\dagger)}^T \bm \Delta \bm P_{\bm V_2}}_F^2 \\
    &= \norm{\bm X^\dagger \bm \Delta^T \bm P_{\bm U_2} \bm \Delta \bm P_{\bm V_2}}_F^2 + \norm{\bm P_{\bm U_2} \bm \Delta \bm P_{\bm V_2} \bm \Delta^T \bm X^\dagger}_F^2 + \norm{\bm P_{\bm U_2} \bm \Delta {(\bm X^\dagger)}^T \bm \Delta \bm P_{\bm V_2}}_F^2 . \numberthis \label{equ:Rg28} 
\end{align*}
\hlnew{Using the definition of $\bm E$ in} (\ref{equ:E}), (\ref{equ:Rg28}) can be represented as 
\begin{align*}
    \norm{\bm R_{2\bm X}(\bm \Delta)}_F^2 &= \norm{\bm U_1 \bm \Sigma_1^{-1} \bm E_{21}^T \bm E_{22} \bm V_2^T}_F^2 + \norm{\bm U_2 \bm E_{22} \bm E_{12}^T \bm \Sigma_1^{-1} \bm V_1^T}_F^2 + \norm{\bm U_2 \bm E_{21} \bm \Sigma_1^{-1} \bm E_{12} \bm V_2^T}_F^2 \\
    &= \norm{\bm \Sigma_1^{-1} \bm E_{21}^T \bm E_{22}}_F^2 + \norm{\bm E_{22} \bm E_{12}^T \bm \Sigma_1^{-1}}_F^2 + \norm{\bm E_{21} \bm \Sigma_1^{-1} \bm E_{12}}_F^2 , \numberthis \label{equ:Rg29} 
\end{align*}
where the second equality stems from Lemma~\ref{lem:semi}. 
Using Lemma~\ref{lem:norm} and the fact that $\norm{\bm \Sigma_1^{-1}}_2 = 1/\sigma_r$, we can bound the RHS of (\ref{equ:Rg29}) by
\begin{align*}
    \norm{\bm R_{2\bm X}(\bm \Delta)}_F^2 &\leq \frac{1}{\sigma_r^2} \Bigl( \norm{\bm E_{21}}_F^2 \norm{\bm E_{22}}_F^2 + \norm{\bm E_{22}}_F^2 \norm{\bm E_{12}}_F^2 + \norm{\bm E_{12}}_F^2 \norm{\bm E_{21}}_F^2 \Bigr) . \numberthis \label{equ:Rg23a}
\end{align*}
\begin{lemma}[\hlnew{Chebyshev's sum inequality} \cite{hardy1952inequalities}] 
\label{lem:abc}
For any $a,b,c \in \R$, we have
\begin{align*}
    3(ab+bc+ca) \leq (a+b+c)^2 .
\end{align*}
\end{lemma}
\noindent Applying Lemma~\ref{lem:abc} to (\ref{equ:Rg23a}) with $a=\norm{\bm E_{21}}_F^2, b=\norm{\bm E_{22}}_F^2$ and $c=\norm{\bm E_{12}}_F^2$, we obtain
\begin{align*}
    \norm{\bm R_{2\bm X}(\bm \Delta)}_F^2 &\leq \frac{1}{\sigma_r^2} \frac{\bigl( \norm{\bm E_{21}}_F^2 + \norm{\bm E_{22}}_F^2 + \norm{\bm E_{12}}_F^2 \bigr)^2}{3} \\
    &\leq \frac{\bigl( \norm{\bm E_{11}}_F^2 + \norm{\bm E_{12}}_F^2 + \norm{\bm E_{21}}_F^2 + \norm{\bm E_{22}}_F^2 \bigr)^2 }{3 \sigma_r^2} = \frac{\norm{\bm E}_F^4}{3 \sigma_r^2} = \frac{\norm{\bm \Delta}_F^4}{3 \sigma_r^2} , \numberthis \label{equ:Rg23}
\end{align*}
where the last equation stems from $\norm{\bm E}_F = \norm{\bm U^T \bm \Delta \bm V}_F = \norm{\bm \Delta}_F$.
From (\ref{equ:Rg23}), \hlnew{taking the square root and then taking the supremum yield}
\begin{align} \label{equ:Rg26}
    \sup_{\norm{\bm \Delta}_F = \epsilon} \norm{\bm R_{2\bm X}(\bm \Delta)}_F &\leq \frac{\norm{\bm \Delta}_F^2}{\sigma_r \sqrt{3}} .
\end{align}
To show that (\ref{equ:Rg26}) implies (\ref{equ:Rg25}), we \hlnew{describe a particular choice of $\bm \Delta$ such that the inequality holds}. Let us choose
\begin{align*}
    \bm \Delta (\epsilon) \triangleq \frac{\epsilon}{\sqrt{3}} \bigl(\bm u_r \bm v_{r+1}^T + \bm u_{r+1} \bm v_r^T + \bm u_{r+1} \bm v_{r+1}^T \bigr) , 
\end{align*}
where $\bm u_r, \bm u_{r+1}, \bm v_r$, and $\bm v_{r+1}$ are the corresponding left and right singular vectors of $\bm X$. Similar to (\ref{equ:c1_1}), one can verify that $\norm{\bm \Delta (\epsilon)}_F = \epsilon$. In addition, from Proposition~\ref{pro:stewart}, we have
\begin{align} \label{equ:Rg30}
    &\bm E_{12} = \bm U_1^T \bm \Delta (\epsilon) \bm V_2 = \frac{\epsilon}{\sqrt{3}} \bm e^r_r (\bm e^{n-r}_{1})^T , ~
    \bm E_{21} = \bm U_2^T \bm \Delta (\epsilon) \bm V_1 = \frac{\epsilon}{\sqrt{3}} \bm e^{m-r}_{1} (\bm e^r_{r})^T , ~
    \bm E_{22} = \bm U_2^T \bm \Delta (\epsilon) \bm V_2 = \frac{\epsilon}{\sqrt{3}} \bm e^{m-r}_{1} (\bm e^{n-r}_{1})^T .
\end{align}
Substituting (\ref{equ:Rg30}) back into (\ref{equ:Rg29}) yields
\begin{align*}
    &\norm{\bm R_{2\bm X} \bigl(\bm \Delta (\epsilon) \bigr)}_F^2 \\
    &\qquad = \norm{\frac{\epsilon^2}{3} \bm \Sigma_1^{-1} \bm e^r_r (\bm e^{m-r}_{1})^T \bm e^{m-r}_{1} (\bm e^{n-r}_{1})^T}_F^2 + \norm{\frac{\epsilon^2}{3} \bm e^{m-r}_{1} (\bm e^{n-r}_{1})^T \bm e^{n-r}_{1} (\bm e^r_r)^T \bm \Sigma_1^{-1}}_F^2 + \norm{\frac{\epsilon^2}{3} \bm e^{m-r}_{1} \bm e_{r}^T \bm \Sigma_1^{-1} \bm e_r (\bm e^{n-r}_{1})^T}_F^2 \\
    &\qquad = \frac{\epsilon^4}{9} \biggl( \frac{1}{\sigma_r^2} + \frac{1}{\sigma_r^2} + \frac{1}{\sigma_r^2} \biggr) = \frac{\norm{\bm \Delta}_F^4}{3\sigma_r^2} . \qquad \qquad \text{(since } \norm{\bm \Delta (\epsilon)}_F = \epsilon \text{)}
\end{align*}
Therefore, the equality in (\ref{equ:Rg26}) holds when $\bm \Delta = \bm \Delta (\epsilon)$, for any $\epsilon>0$. This completes our proof of the theorem.

\subsection{Proof of Lemma~\ref{lem:max}}

Since $f(\bm x)-g(\bm x) \leq \abs{f(\bm x)-g(\bm x)}$, we have $f(\bm x) \leq \abs{f(\bm x)-g(\bm x)} + g(\bm x)$.
Taking the supremum yields
\begin{align*}
    \sup_{\bm x \in \C} f(\bm x) &\leq \sup_{\bm x \in \C} \biggl\{ \abs{f(\bm x)-g(\bm x)} + g(\bm x) \biggr\} \\
    &\leq \sup_{\bm x \in \C} \abs{f(\bm x)-g(\bm x)} + \sup_{\bm x \in \C} g(\bm x).
\end{align*}
Thus, \hlnew{we have}
\begin{align} \label{equ:sup1}
    \sup_{\bm x \in \C} f(\bm x) - \sup_{\bm x \in \C} g(\bm x) \leq \sup_{\bm x \in \C} \abs{f(\bm x) - g(\bm x)} . 
\end{align}
Changing the roles of $f$ and $g$, we also obtain
\begin{align} \label{equ:sup2}
    \sup_{\bm x \in \C} g(\bm x) - \sup_{\bm x \in \C} f(\bm x) \leq \sup_{\bm x \in \C} \abs{f(\bm x) - g(\bm x)} . 
\end{align}
Our inequality follows on combining (\ref{equ:sup1}) and (\ref{equ:sup2}).


    
    
    



\bibliographystyle{model1-num-names}

\bibliography{fulltitles,cas-refs}

\end{document}